\documentclass[12pt]{article}
\usepackage[utf8]{inputenc}
\usepackage{graphicx} 
\usepackage{amsfonts}
\usepackage{amsmath}
\usepackage{amssymb}
\usepackage{amsthm}
\usepackage{tikz}
\usetikzlibrary{positioning,arrows.meta,cd}
\usetikzlibrary{angles,quotes}
\usepackage{verbatim}
\usepackage{booktabs}
\usepackage[english]{babel}
\usepackage{mathtools}
\usepackage{authblk}

\theoremstyle{definition}
\newtheorem{theorem}{Theorem}[section]
\newtheorem{proposition}[theorem]{Proposition}
\newtheorem{lemma}[theorem]{Lemma}

\newtheorem{definition}[theorem]{Definition}

\title{The Complexity of $L(3n,1)$}
\author{William Hobkirk}
\author{Jonathan Spreer}
\affil{School of Mathematics and Statistics, University of Sydney NSW 2006, Australia, \texttt{william.hobkirk@sydney.edu.au}, \texttt{jonathan.spreer@sydney.edu.au}}
\date{}

\begin{document}

\maketitle

\begin{abstract}
    We show that any triangulation of a lens space $L(3n,1)$, $n \geq 2$, must have 
    at least $3n-3$ tetrahedra. This lower bound coincides with the number of tetrahedra 
    in the corresponding minimal layered triangulations -- confirming the conjecture by 
    Jaco-Rubinstein and Matveev that they are minimal in these cases. 
    For even $n$, this bound also follows from work by Jaco, Rubinstein and Tillmann 
    proving the analogous bound for the lens spaces $L(2n,1)$.   
    We extend their method of using the Euler characteristic of embedded surfaces 
    representing $\mathbb{Z}_2$-homology classes, by applying the theory of folded 
    surfaces describing $\mathbb{Z}_3$-homology classes.
\end{abstract} 



\section{Introduction}
\label{ssec:intro}

The triangulation complexity $c(M)$ of a closed $3$-manifold $M$ is the smallest number of tetrahedra in any (generalised) triangulation of $M$. It is a useful organising principle, analogous to the crossing number of a knot: for every $k>0$, only finitely many closed topological $3$-manifolds have triangulation complexity at most $k$. In principle, these manifolds can be enumerated, and censuses have been produced using computational tools such as \emph{Regina} \cite{regina}. For example, there are exactly $13\,399$ closed, orientable, irreducible $3$-manifolds of complexity at most $11$ \cite{Matveev2005Tabulation}. More recent work has extended the tabulation of closed, orientable, irreducible $3$-manifolds to complexity $13$~\cite{MatveevTarkaev2020}. 

Despite these computations, determining the triangulation complexity of a given 
$3$-manifold remains difficult. Census methods yield exact answers only for finitely 
many $3$-manifolds, and there is typically a substantial gap between existing lower 
bounds and the number of tetrahedra in the smallest known triangulations. 

Lower bounds on triangulation complexity have been developed using several different topological and geometric invariants. For a hyperbolic $3$-manifold, its hyperbolic volume gives the lower bound $c(M)\geq \operatorname{Vol}(M)/v_3$, where $v_3$ is the volume of the regular ideal tetrahedron, see, for instance, \cite{Thurston1979}. Volume estimates can be used to derive lower bounds for particular hyperbolic families, including closed hyperbolic manifolds and $2$-bridge link complements \cite{MatveevPetronioVesnin2009,MorganSpreer2026}. Cha obtained lower bounds from Cheeger--Gromov $\rho$-invariants, including examples for which volume-based bounds are arbitrarily weaker \cite{Cha2016Surgery}. Lackenby and Purcell related the complexity of hyperbolic surface bundles to translation lengths of their monodromies \cite{LackenbyPurcell2024Fibred}, and established bounds up to universal multiplicative constants for all elliptic and Sol $3$-manifolds, including lens spaces \cite{lackenby2019}. Further bounds, in some cases matching explicit upper bounds up to small additive constants, have been obtained for infinite families of Dehn fillings of cusped $3$-manifolds \cite{JacoRubinsteinSpreerTillmann2025}. 

A particularly successful method for obtaining lower bounds depends on the 
$\mathbb Z_2$-Thurston norm. Jaco, Rubinstein and Tillmann related the complexity 
of a triangulation to the topology of normal surfaces representing a non-trivial 
$\mathbb Z_2$-homology class \cite{JacoRubinsteinTillmann2013}. They used this method 
to determine the {\em exact} triangulation complexity for infinite families of 
$3$-manifolds, including lens spaces $L(2n,1)$, $n\geq 2$ \cite{Jaco_2009}, and other 
spherical $3$-manifolds \cite{Jaco_2011}. More recently, this theory was developed further, 
yielding exact results for additional infinite families of $3$-manifolds \cite{Jaco_2020}.

Despite all this progress, significant questions remain open: In 
\cite{jaco2006layeredtriangulations}, Jaco and Rubinstein constructed
\emph{minimal layered triangulations} for all lens spaces, providing explicit upper bounds 
for their triangulation complexity. They, as well as Matveev, conjectured that all of these 
minimal layered triangulations realise the triangulation complexity for their underlying
manifolds, see \cite[§6, p. 54]{jaco2006layeredtriangulations} and 
\cite[§2.3.3, item E, p. 77]{Matveev2007}. To date, this conjecture has 
only been proven for a small number of one-parameter infinite families \cite{Jaco_2009, Jaco_2011, Jaco_2020}.

The purpose of this article is to make further progress on this conjecture, by applying an 
additional approach to proving triangulation complexity results. In particular, we extend 
the homological approach, introduced by Jaco, Rubinstein and Tillmann in \cite{Jaco_2009},
beyond $\mathbb Z_2$. Motivated by Turaev's study of folded surfaces representing torsion 
homology classes \cite{Turaevfnhom3mfld2007}, we use $\mathbb Z_3$-folded surfaces and their 
associated norm, developed in \cite{hobkirk2026}, to obtain a sharp complexity bound for 
$L(3n,1)$ for every $n\geq 2$. The resulting bound agrees with the number of tetrahedra in 
the minimal layered triangulation. In particular, we prove the following result:

\begin{theorem}
    \label{thm:main}
    The triangulation complexity of the lens space $L(3n,1)$ is $3n-3$ for all integers $n\geq 2$.
\end{theorem}

This suggests a route towards complexity bounds using $\mathbb Z_p$-torsion homology for
other odd $p$, not just for lens spaces, but also for other $3$-manifolds, reducing the 
topological restrictions on homological lower-bound methods.

\paragraph{Structure of this paper} In the remainder of this introduction we give some background on the origin of this research question. Section~\ref{ssec:foldedsurfaces} provides the necessary background for the proof of Theorem~\ref{thm:main} and summarises the theory of folded surfaces developed in \cite{hobkirk2026}. Section~\ref{sec:proofsec} contains the proof of Theorem~\ref{thm:main}, before we discuss future work in Section~\ref{sec:discussion}.

%
%


%
%

\paragraph{The origins of this article}

This section is inspired by the blog post ``Let the Diners Into the Kitchen'', written by Noah Giansiracusa on Terry Tao's blog \cite{giansiracusa2026diners}. It serves the purpose of giving this research story a human face.

The idea of using $\mathbb{Z}_p$-homology classes, $p > 2$, to determine triangulation complexity of $3$-manifolds dates back to discussions between Hyam Rubinstein, Stephan Tillmann and the second author in 2019, shortly after the second author moved to Sydney.

Back then, we were not sure how to reconcile having a singular locus in a class representative with a non-orientable non-singular surface representing $\mathbb{Z}_2$-classes. We had the initial intuition of interpreting a non-orientable surface as a bounded planar surface with antipodal points on its boundary component identified, but failed to capitalise on this viewpoint. 

We ultimately abandoned the project, until the second author -- encouraged by discussions with Daniele Celoria -- started thinking about the project again in 2022, now with the help of the first author who had just started his PhD and was brave enough to take on the literature on folded surfaces.

Under the guidance of the first author, the project then pivoted towards a more thorough investigation into a simple, geometric representation of torsion homology classes, merely guided in the background by the initial motivation to prove triangulation complexity bounds. This work ultimately resulted in a comprehensive body of theory, published in \cite{hobkirk2026}, that discusses representations of torsion classes in a ``triangulation-friendly'' way. This theory uncovered many interesting facts about lens spaces in surprisingly simple and vivid ways, as is demonstrated on the first author's website \cite{hobkirkWebsite}. 

Moreover, it turns out that the number of crosscaps in the $\mathbb{Z}_2$ case, as well as the number of singular locus components in the $\mathbb{Z}_p$, $p>2$, case is less important: in both cases, adding handles to the orientable part of the surface can be exchanged for the number of crosscaps or components. This observation connects our theory to the existing theory of folded surfaces by Turaev \cite{Turaevfnhom3mfld2007}.

Ongoing discussions and presentations of early versions of this work in the years 2025 and 2026 have further refined our understanding of what we actually want to present, and which parts of the theory and its applications we will leave for future work, or work by others. Most significantly, we decided to present the theory necessary to prove the norm of $\mathbb{Z}_3$-classes in lens spaces in a separate publication. This was partly to give the first author the appropriate level of credit for his work, and partly because this modular presentation follows a very natural separation of the two parts of a triangulation complexity proof: a description of a topological quantity in a family of 3-manifolds, and a combinatorial argument translating this quantity into a lower bound on the number of tetrahedra in any triangulation of these manifolds. 

\subsection*{Acknowledgements and AI declaration}

Work of the second author is supported by the Australian Research Council under the Discovery Project scheme, grant no. DP220102588. The authors would like to thank Daniele Celoria and Stephan Tillmann for helpful discussions. 

We used AI for proofreading of this article. We also used AI to organise the combinatorics in Proposition~\ref{prop:chargeineqequiv}. We completed the remainder of the proof and the rest of the article without using AI.

\section{Triangulations and $\mathbb{Z}_3$-folded surfaces}
\label{ssec:foldedsurfaces}

This article is about lens spaces. This is a family of $3$-manifolds that is obtained by gluing the boundaries of two solid tori together. Fix a meridian $\mu$ and a longitude $\lambda$ of the second solid torus. The lens space $L(p,q)$, $p$ and $q$ co-prime, is obtained by the identification that glues a meridian of the first solid torus to a curve of slope $p\lambda + q\mu$ on the boundary of the second solid torus. The homeomorphism type of the resulting closed $3$-manifold does not depend on the particular choices of meridian and longitude.

We describe lens spaces, and more generally, $3$-manifolds, using what is called a {\em generalised triangulation}, sometimes also called a {\em pseudo-simplicial triangulation}. Starting with a collection of disjoint tetrahedra $\{ \Delta_1 , \ldots , \Delta_n\}$, and gluings between pairs of their boundary triangles $\Phi_i$, $ 1\leq i \leq 2n$, we obtain a space $\mathcal{T} = \{ \Delta_1 , \ldots , \Delta_n \} / \Phi_i$ using the quotient topology. As a result of the triangle gluings, vertices, edges and triangles become identified, and we refer to an equivalence class under this identification as a vertex, edge or triangle of the triangulation $\mathcal{T}$, and similarly so with tetrahedra. Given a triangulation $\mathcal{T}$, we denote its set of $i$-dimensional faces by $\mathcal{T}^{(i)}$, $0 \leq i \leq 3$. Given an edge $e \in \mathcal{T}^{(1)}$, the degree $\deg (e)$ of $e$ is the number of edges of $\{ \Delta_1 , \ldots , \Delta_n\} $ that become identified under the triangle gluings to form $e$. 

We say that a triangulation $\mathcal{T}$ is {\em valid} if no edge becomes identified with itself in reverse as a result of the triangle gluings. Given a valid triangulation $\mathcal{T}$, the {\em link} of a vertex $v \in \mathcal{T}^{(0)}$ is the boundary of a small neighbourhood around $v$. We say that $\mathcal{T}$ is a {\em triangulated $3$-manifold} if all of its vertex links are triangulated $2$-dimensional spheres. 

Given a triangulated $3$-manifold $\mathcal{T}$, we denote its underlying set by $|\mathcal{T}|$. Given a $3$-manifold $M$, we say that $\mathcal{T}$ {\em triangulates} $M$ if $|\mathcal{T}| \cong M$.
Readers accustomed to {\em combinatorial manifolds} may find it helpful to note that the second derived subdivision of a generalised triangulation of a $3$-manifold is always a {\em combinatorial manifold}, as defined in \cite{Spreer2011}.

It is a result of Moise that every closed $3$-manifold admits a (generalised) triangulation \cite{moise:1952}. This motivates the following natural question: what is the \emph{simplest} triangulation of a given $3$-manifold $M$? The standard measure of the complexity of a triangulation is the number of tetrahedra in it. Following this, the \emph{triangulation complexity} $c(M)$ of a $3$-manifold $M$ is the minimum number of tetrahedra in any triangulation of $M$. From this perspective the simplest triangulations of $M$ are those with exactly $c(M)$ tetrahedra, which we refer to as \emph{minimal}. Note that the notion of triangulation complexity for orientable, irreducible $3$-manifolds coincides with Matveev's notion of complexity except for the $3$-sphere, real projective $3$-space, and $L(3,1)$ \cite{Matveev2007}. 

\medskip

For the remainder of this section, we recall definitions and results from~\cite{hobkirk2026} 
that we rely on to prove Theorem~\ref{thm:main}. The reader is encouraged to consult this 
source for a more precise and general treatment.

We start by defining $\mathbb{Z}_3$-folded surfaces constructively in Definition~\ref{dfn:z3foldedsurfaces}. See \cite[Definition 3.1]{hobkirk2026} for a more intrinsic characterisation.

\begin{definition}[$\mathbb{Z}_3$-folded surfaces]\label{dfn:z3foldedsurfaces}
    Let $X_{cut}$ be any compact oriented $2$-manifold with boundary. Let $q:\partial X_{cut}\rightarrow \text{Sing}(X)$ be a $3$-sheeted covering map of a closed oriented $1$-manifold $\text{Sing}(X)$ with the property that $q$ is orientation preserving with respect to the orientation of $\partial X_{cut}$ induced by the orientation of $X_{cut}$. A \emph{$\mathbb{Z}_3$-folded surface} is the space $X$ given by the topological quotient of $X_{cut}$ by the equivalence relation on $\partial X_{cut}$ induced by $q$. The \emph{singular locus} of $X$ is $\text{Sing}(X)$. 
\end{definition}

In Definition~\ref{dfn:fundamentalclass} we describe the \textit{fundamental class} of a $\mathbb{Z}_3$-folded surface which is central to the motivating application of these objects for studying $3$-manifolds $M$ with non-trivial $H_2(M;\mathbb{Z}_3)$. To do this we first introduce some notation. Commutative diagram~\eqref{comm:diag} gives labels for the inclusion, quotient and connecting homomorphisms of the long exact sequences of homology of the pairs $(X_{cut},\partial X_{cut})$ and $(X,\text{Sing}(X))$ as well as the long exact sequence of homology induced by the coefficient short exact sequence of Equation~\eqref{eq:shortcoeffseq}. The commutativity of the squares in~\eqref{comm:diag} follows from the naturality of these sequences.

\begin{equation}\label{eq:shortcoeffseq}
    0 \rightarrow \mathbb{Z} \xrightarrow{\cdot 3}\mathbb{Z} \rightarrow \mathbb{Z}/3\mathbb{Z} \rightarrow 0
\end{equation}

\begin{equation}\label{comm:diag}
    \begin{tikzcd}[column sep=large, inner sep=large]
        & H_2(X_{cut},\partial X_{cut};\mathbb{Z})
          \arrow[r,"d_{(X_{cut},\partial X_{cut})}"]
          \arrow[d,"q_*"]
        & H_1(\partial X_{cut};\mathbb{Z})
          \arrow[d,"q_*"] \\
        & H_2(X,\text{Sing}(X);\mathbb{Z})
          \arrow[r,"d_{(X,\text{Sing}(X))}"]
          \arrow[d,"j_{\mathbb{Z}_3}"]
        & H_1(\text{Sing}(X);\mathbb{Z})
          \arrow[d,"j_{\mathbb{Z}_3}"] \\
        H_2(X;\mathbb{Z}_3)
          \arrow[r,"j_{(X,\text{Sing}(X))}"]
        & H_2(X,\text{Sing}(X);\mathbb{Z}_3)
          \arrow[r,"d_{(X,\text{Sing}(X))}"]
        & H_1(\text{Sing}(X);\mathbb{Z}_3)
    \end{tikzcd}
\end{equation}

\begin{definition}[The fundamental class of a $\mathbb{Z}_3$-folded surface]\label{dfn:fundamentalclass}
    Let $X$ be a $\mathbb{Z}_3$-folded surface as in Definition~\ref{dfn:z3foldedsurfaces}. 
    Let $[X_{cut}]\in H_2(X_{cut},\partial X_{cut};\mathbb{Z})$, 
    $[\partial X_{cut}] \in H_1(\partial X_{cut};\mathbb{Z})$ and 
    $[\text{Sing}(X)]\in H_1(\text{Sing}(X);\mathbb{Z})$ denote the fundamental classes of 
    $X_{cut}, \partial X_{cut}$, and $\text{Sing}(X)$, respectively. By definition, 
    $d_{(X_{cut},\partial X_{cut})}([X_{cut}]) = [\partial X_{cut}]$ and 
    $q_*([\partial X_{cut}]) = 3[\text{Sing}(X)]$. Since $q_*([\partial X_{cut}])$ 
    is in the kernel of 
    $j_{\mathbb{Z}_3}:H_1(\text{Sing}(X);\mathbb{Z}) \to H_1(\text{Sing}(X);\mathbb{Z}_3)$ it follows that $j_{\mathbb{Z}_3} \circ q_*([X_{cut}])$ is in the kernel of 
    $d_{(X,\text{Sing}(X))}:H_2(X,\text{Sing}(X);\mathbb{Z}_3)\rightarrow H_1(\text{Sing}(X);\mathbb{Z}_3)$. 
    Hence by exactness of the sequence of the pair $(X,\text{Sing}(X))$ there must exist some 
    class $[X]\in H_2(X;\mathbb{Z}_3)$ such that 
    $j_{(X,\text{Sing}(X))}([X]) = j_{\mathbb{Z}_3} \circ q_*([X_{cut}])$. Moreover, since 
    $H_2(\text{Sing}(X);\mathbb{Z}_3) = 0$ we can deduce that the map 
    $j_{(X,\text{Sing}(X))}:H_2(X;\mathbb{Z}_3) \rightarrow H_2(X,\text{Sing}(X);\mathbb{Z}_3)$ 
    is injective and hence that $[X]$ is the unique class with this property.    
    We define $[X]$ to be the \emph{fundamental class of $X$}. 
\end{definition}

If $X$ is embedded in a closed orientable $3$-manifold $M$ then we say that $X$ represents a class $x\in H_2(M;\mathbb{Z}_3)$ if $x$ is the image of the fundamental class of $X$ under the inclusion map. Definition~\ref{dfn:canonicalz3folded} recalls the construction of a \textit{canonical $\mathbb{Z}_3$-folded surface} representing the dual of some class $\varphi\in H^1(M;\mathbb{Z}_3)$ and which depends only on a one-vertex triangulation $\mathcal{T}$ of $M$. This proves the existence of $\mathbb{Z}_3$-folded surfaces representing every class of $H_2(M;\mathbb{Z}_3)$. More details on this construction and its generalisation to $\mathbb{Z}_p$-folded surfaces for any odd $p$ can be found in~\cite[§5, p. 39]{hobkirk2026}.

\begin{definition}[Canonical $\mathbb{Z}_3$-folded surfaces]\label{dfn:canonicalz3folded}

    Let $\mathcal{T}$ be a one-vertex triangulation of a closed, oriented $3$-manifold $M$. Let $\varphi \in H^1(M;\mathbb{Z}_3)$.
    The canonical $\mathbb{Z}_3$-folded surface $S_{\varphi}(\mathcal{T})$ is a $\mathbb{Z}_3$-folded surface that is constructed as follows.
    
    First note that each edge of this triangulation is a loop in $M$ since both endpoints 
    coincide at the single vertex. Via the universal coefficient theorem we identify 
    $H^1(M;\mathbb{Z}_3) = \text{Hom}(H_1(M;\mathbb{Z}), \mathbb{Z}_3)$. Choose an orientation 
    for each edge $e$ such that $\varphi([e]) \in \{0,1\} \subset \mathbb{Z}/3\mathbb{Z}$. We 
    refer to edges $e$ such that $\varphi([e]) = 0$ as \textit{even} edges and edges such that
    $\varphi([e]) = 1$ as \textit{odd edges}. Since the orientation of an even edge is 
    arbitrary we do not indicate it in our figures. The relations between the 
    $H_1(M;\mathbb{Z})$ classes of edges on the boundary of each triangle imply that there are
    precisely four ways that $\varphi$ can colour the boundary of a triangle and five ways 
    $\varphi$ can colour the boundary of a tetrahedron in $\mathcal{T}$ up to symmetry. 
    These colourings are given in Figure~\ref{fig:admissiblecolourings}.
    
    Partition the tetrahedra of $\mathcal{T}$ into types $1$, $2$, $3$ and $4$ which have $2,3,6$ and $1$ pre-images of even edges on their boundaries, respectively. We construct $S_{\varphi}(\mathcal{T})$ by specifying its intersection with each tetrahedron as in Figure~\ref{fig:elementarypieces} such that the elementary pieces intersect pre-images of odd edges in a single point and do not intersect pre-images of even edges. A tetrahedron of type $1$ intersects $S_{\varphi}(\mathcal{T})$ in a quadrilateral. The intersection with a type $2$ tetrahedron is a triangle. The intersection with a type $3$ tetrahedron is empty. Finally, the intersection with a type $4$ tetrahedron is $K_{3,1} \times [0,1]$ where $K_{3,1}$ is the tripod graph. Orient these elementary pieces such that the pairing with the orientation of the odd edges at each point of intersection is positive with respect to the orientation of $M$. These elementary pieces can be specified in such a way that the face pairing maps restrict to orientation reversing maps on the boundary of each elementary piece. It follows that, in $M$, these elementary pieces glue up into a $\mathbb{Z}_3$-folded surface $S_{\varphi}(\mathcal{T})$ and the signed count of intersections with each edge guarantees that $S_{\varphi}(\mathcal{T})$ represents the Poincar\'e dual of $\varphi$.
    $S_{\varphi}(\mathcal{T})$ is unique up to ambient isotopy.

\begin{figure}
    \centering
    \includegraphics[width=0.5\linewidth]{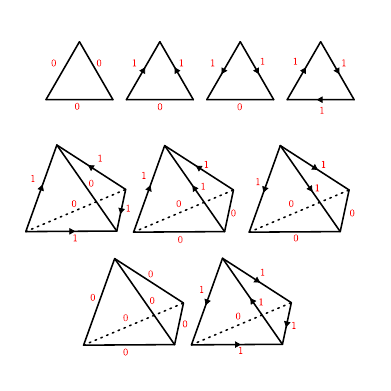}
    \caption{Admissible colourings of triangles and tetrahedra up to symmetry.}
    \label{fig:admissiblecolourings}
\end{figure}

\begin{figure}
    \centering
    \includegraphics[width=0.5\linewidth]{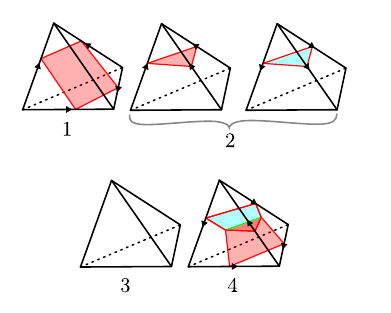}
    \caption{Elementary pieces of a canonical $\mathbb{Z}_3$-folded surface in tetrahedra of types $1,2,3$ and $4$.}
    \label{fig:elementarypieces}
\end{figure}

\end{definition}

Lemma~\ref{lem:nonseparating} is a specialisation of \cite[Lemma 5.2, p. 44]{hobkirk2026} that is of particular relevance to our paper, proving that canonical $\mathbb{Z}_3$-folded surfaces have the special property that they do not separate the ambient $3$-manifold. We call $\mathbb{Z}_3$-folded surfaces that have this property \emph{non-separating}.

\begin{lemma}\label{lem:nonseparating}
    For any closed, orientable, connected $3$-manifold $M$, any one-vertex triangulation $\mathcal{T}$ of $M$, and $\varphi \in H^1(M;\mathbb{Z}_3)$ non-trivial, the canonical $\mathbb{Z}_3$-folded surface $S_{\varphi}(\mathcal{T})$ is non-separating.
\end{lemma}
    
\begin{proof}
    By checking each tetrahedron type in Figure~\ref{fig:elementarypieces} we can verify that every point in the complement of each elementary piece is path connected to at least one vertex pre-image. Passing to the quotient and noting that there is a single vertex class in $M$ this fact is sufficient to show that $S_{\varphi}(\mathcal{T})$ does not separate the ambient space. 
\end{proof}

Non-separating $\mathbb{Z}_3$-folded surfaces were studied extensively in~\cite{hobkirk2026}, 
where it was shown that their Euler characteristic cannot be arbitrarily large. 
Definition~\ref{dfn:chi3} reproduces the $\chi_3$ function which characterises the 
non-separating representatives of a fixed second $\mathbb{Z}_3$-homology class.

\begin{definition}\label{dfn:chi3}
    Let $M$ be a closed orientable $3$-manifold and let $x \in H_2(M;\mathbb{Z}_3)$. Define $\chi_3(x) \in \mathbb{Z}$ to be the maximum Euler characteristic of any $\mathbb{Z}_3$-folded surface $X \subset M$ representing $x$ with the property that $M\setminus X$ is connected.
\end{definition}

The main result of~\cite{hobkirk2026} is Theorem~\ref{thm:maxchinonsep} which determines the maximum Euler characteristic of a non-separating non-empty $\mathbb{Z}_3$-folded surface that embeds in $L(3n,1)$ by generalising techniques developed by Bredon and Wood in~\cite{bredonwood1969} in the context of non-orientable surfaces. 

\begin{theorem}[Theorem~4.21 in \cite{hobkirk2026}]\label{thm:maxchinonsep}
    Let $x$ be a non-trivial class in $H_2(L(3n,1);\mathbb{Z}_3)$, $n \geq 1$, then $\chi_3(x) = 2-n$.
\end{theorem}

A corollary of Theorem~\ref{thm:maxchinonsep} and Lemma~\ref{lem:nonseparating} is that if
$\varphi \in H^1(L(3n,1);\mathbb{Z}_3)$ is non-trivial, then 
$\chi(S_{\varphi}(\mathcal{T})) \leq 2-n$ for all one-vertex triangulations $\mathcal{T}$. 
Section~\ref{sec:proofsec} uses this fact to determine the triangulation complexity 
of $L(3n,1)$.

\section{Proof of Theorem~\ref{thm:main}}\label{sec:proofsec}
\subsection{Overview}

Let $n \geq 2$ be an integer. The \emph{minimal layered triangulation} of $L(3n,1)$ defined in \cite{jaco2006layeredtriangulations} has exactly $3n-3$ tetrahedra and so by definition $c(L(3n,1)) \leq 3n-3$ \cite[§2.5, p. 161]{Jaco_2009}. Therefore it is sufficient to prove the lower bound $c(L(3n,1))\geq 3n-3$ for Theorem~\ref{thm:main}. The broad strokes are as follows: In Proposition~\ref{prop:chargeineqequiv} we apply standard arguments from triangulation combinatorics to manipulate the desired lower bound on the number of tetrahedra into an equivalent inequality that depends only on the combinatorics and topology of the non-trivial canonical $\mathbb{Z}_3$-folded surface of a minimal triangulation. Lemma~\ref{eq:allAdeg5} determines a lower bound on the gap between the Euler characteristic of this canonical $\mathbb{Z}_3$-folded surface and the maximum given by Theorem~\ref{thm:maxchinonsep} using disc compressions. Lemma~\ref{lem:psiineq} applies this to derive a sufficient condition for the inequality of Proposition~\ref{prop:chargeineqequiv} and hence the desired lower bound. Lemma~\ref{lem:psiineq} depends on a \emph{surplus} function which assigns a rational number to each even edge in the minimal triangulation depending on the number and types of tetrahedra it is contained in, counted by multiplicity. The sufficient condition for a proof of Theorem~\ref{thm:main} is that the \emph{total surplus} (the sum of the surplus function over all the even edges in the triangulation) is non-negative. Lemma~\ref{lem:classifydeficitedges} uses the classification of edges of low degree from \cite{Jaco_2009} to eliminate all but a very small number of cases in which an even edge may have negative surplus (or \emph{deficit}). We then compute a lower bound on the global surplus in terms of combinatorial properties of a subcomplex of the triangulation in Lemma~\ref{lem:sufficientconditionpinchedhandle}. This allows us to decompose the global surplus inequality into sufficient conditions that depend on the combinatorics of certain sub-triangulation blocks. Proposition~\ref{prop:desiredbound} establishes the required condition for these blocks. 

\subsection{The surplus inequality}

Let $\mathcal{T}$ be a minimal triangulation of the lens space $L(3n,1)$. It follows from~\cite[Theorem 6.1]{jaco20020efficienttriangulations3manifolds} that $\mathcal{T}$ has one vertex and is $0$-efficient. Let $\varphi$ be either non-trivial class in $H^1(L(3n,1);\mathbb{Z}/3\mathbb{Z})= \mathbb{Z}_3$. Let $X = S_{\varphi}(\mathcal{T})$ be the canonical $\mathbb{Z}_3$-folded surface dual to $\varphi$. Define the following quantities:

\begin{itemize}
    \item Let $\mathfrak{e}$ be the number of even edges of $\varphi$ in $\mathcal{T}$.
    \item Let $A,B,C$ and $D$ be the number of type $1,2,3$ and $4$ tetrahedra in the triangulation, respectively.
    \item For each edge $e\in \mathcal{T}^{(1)}$ let $a_e, b_e, c_e$ and $d_e$ be the number of pre-images of $e$ in tetrahedra of types $1,2,3$ and $4$, respectively.
    \item Let $s(e)= \frac{a_e}{2} + \frac{5b_e}{6} + \frac{2c_e}{3} + d_e$ for each edge $e$. 
\end{itemize}

Proposition~\ref{prop:chargeineqequiv} expresses the conjectured lower bound on the number of tetrahedra in $\mathcal{T}$ in terms of the sum of $s(e)$ over all even edges $e$.

\begin{proposition} \label{prop:chargeineqequiv}
    Let $\mathcal{T}$ be a one-vertex triangulation of $L(3n,1)$.
    The following are equivalent:
    \begin{align}
        \sum_{\varphi(e)=0} s(e) &\geq 3\mathfrak{e} + 3\chi(X) - 3\chi_3([X]) \label{eq:chargeineq} \\
        \#\mathcal{T}^{(3)} &\geq 3n -3 \notag
    \end{align}
\end{proposition}
\begin{proof}
    First note that by counting the number of pre-images of even edges in the boundary of each type of tetrahedron we obtain the following relation:

    \begin{equation}\label{eq:chargeeq}
        \sum_{\varphi(e)=0}s(e) = A + \frac{5}{2}B + 4C + D = |\mathcal{T}^{(3)}| + \frac{3}{2}B + 3C 
    \end{equation}

    It follows from elementary topological arguments that the canonical $\mathbb{Z}_3$-folded surface $X$ has the same Euler characteristic as its complement in $L(3n,1)$. We can also construct a homotopy of $L(3n,1)\setminus X$ onto the subcomplex of $\mathcal{T}$ spanned by the set of even edges by piecing together homotopies in each tetrahedron type that respect the face pairings. The subcomplex spanned by the set of even edges has one vertex, $\mathfrak{e}$ edges, $\frac{B+4C}{2}$ triangles and $C$ tetrahedra. Together, these observations imply

    \begin{equation}\label{eq:eulercharcomp}
        \chi(X) = 1 - \mathfrak{e} + B/2 + C
    \end{equation}

    Rearranging Equation~\eqref{eq:eulercharcomp} gives $\frac{3}{2}B + 3C = 3\mathfrak{e} - 3 + 3\chi(X)$ which we substitute into Equation~\eqref{eq:chargeeq} to obtain
    
    \[\sum_{\varphi(e)=0} s(e) = |\mathcal{T}^{(3)}| + 3\mathfrak{e} - 3 + 3\chi(X)\]
    Hence,
    \[ |\mathcal{T}^{(3)}| = \sum_{\varphi(e)=0} s(e) - 3\mathfrak{e} + 3 - 3\chi(X)\]
    Applying Inequality~\eqref{eq:chargeineq} gives
    
    \[ |\mathcal{T}^{(3)}| \geq 3 - 3\chi_3([X]).\]
    
    Theorem~\ref{thm:maxchinonsep} implies that $\chi_3([X]) = 2-n$ since $[X]$ is non-trivial, and so this is equivalent to the required inequality
    \[|\mathcal{T}^{(3)}| \geq 3 - 6 + 3n = 3n-3\]
\end{proof}

Assume that $\mathcal{T}$ is a minimal triangulation. We seek to prove that Inequality~\eqref{eq:chargeineq} holds on any minimal triangulation $\mathcal{T}$ of $L(3n,1)$. First note that by definition $\chi(X) \leq \chi_3([X])$. We start by proving Lemma~\ref{eq:allAdeg5}.

\begin{lemma}\label{eq:allAdeg5}
    Let $\mathcal{A}$ be the set of even edges of $\mathcal{T}$ with every pre-image in a type $1$ tetrahedron. The following lower bound holds:

    \[\chi_3([X]) - \chi(X) \geq | \mathcal{A} |.\]
\end{lemma}

\begin{proof}
    Let $e \in \mathcal{A}$ and hence $a_e = \deg(e)$. Every pre-image of $e$ must occur in a distinct tetrahedron, otherwise one can find a M\"obius band embedded in $\text{Int}(X)$ which would contradict the fact that $\text{Int}(X)$ is oriented. Hence, the pairings of triangles in which $e$ lies form a triangulation of a $3$-ball and there exists a compression disc of $\text{Int}(X)$ embedded in the interior of this $3$-ball which intersects $e$ transversely. Figure~\ref{fig:compressiondiskA} gives a visualisation of two such compression discs.

    \begin{figure}
        \centering
        \includegraphics[width=0.5\linewidth]{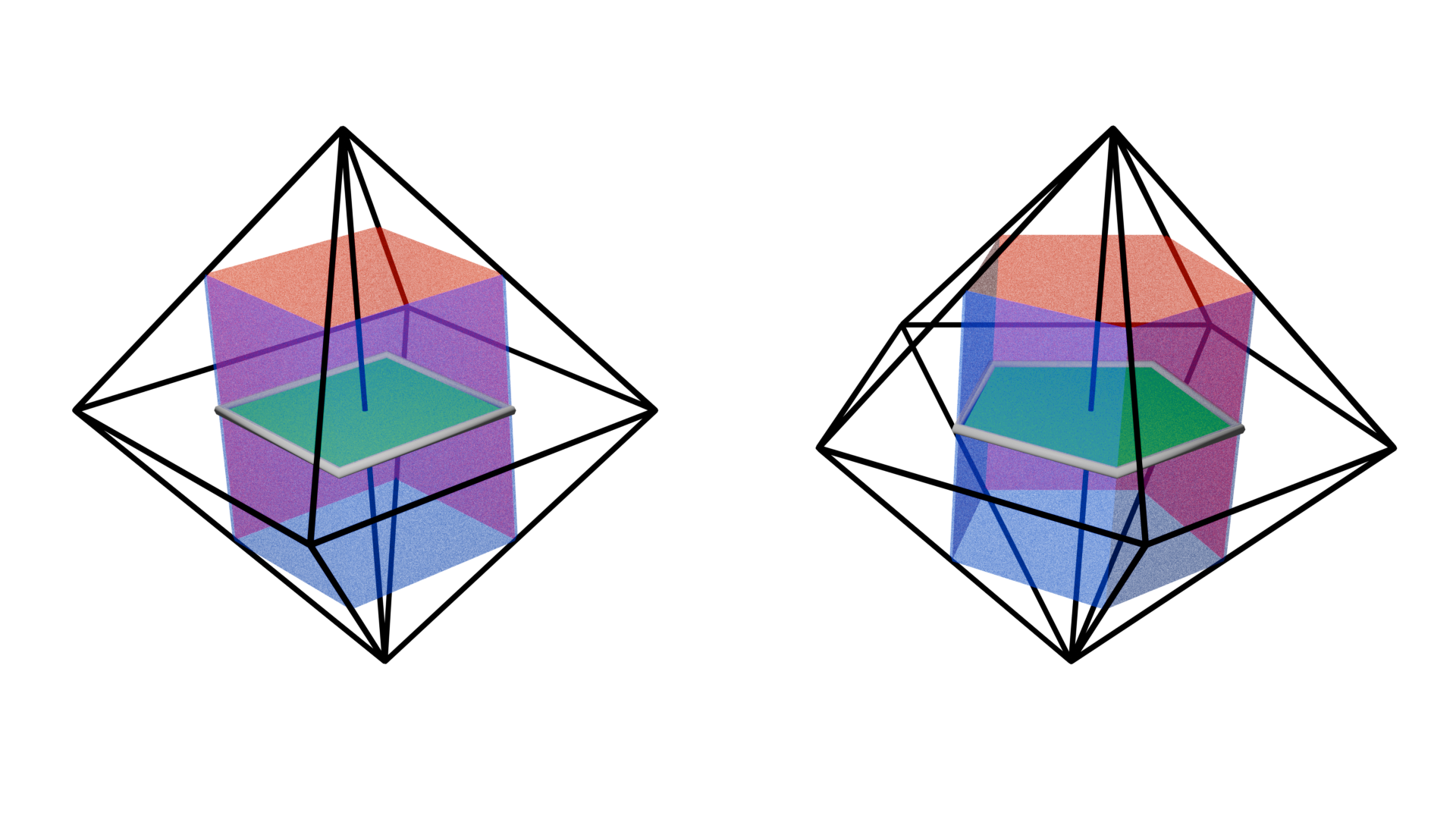}
        \caption{A compression disc of $\text{Int}(X)$ in $3$-ball models of edges in $\mathcal{A}$ of degree $4$ and $5$. These two edges belong to the same subset in the partition of $\mathcal{A}$ since the compression discs lie on the same side of $\text{Int}(X)$.}
        \label{fig:compressiondiskA}
    \end{figure}

    Moreover, we can choose an arbitrary co-orientation of $\text{Int}(X)$ which induces a partition of $\mathcal{A}$ into subsets $\mathcal{A}_+$ (resp. $\mathcal{A}_-$) consisting of the edges in $\mathcal{A}$ whose transverse compression discs lie on the positive (resp. negative) side of $\text{Int}(X)$. Without loss, assume that $|\mathcal{A}_+| \geq | \mathcal{A}_-|$, and hence $|\mathcal{A}_+| \geq |\mathcal{A}|/2$. 

    This choice of co-orientation forces the tetrahedra in the $3$-ball models around each edge
    in $\mathcal{A}_+$ to be pairwise disjoint and so compression discs corresponding to these 
    edges are all disjoint. Performing compression surgeries along these discs does not affect 
    the non-separating property of $X$ since there is at least one pre-image of the single 
    vertex of $\mathcal{T}$ in each complementary region of each $3$-ball model after the 
    surgery as in Figure~\ref{fig:compressiondiskA}.
   
    Hence we can perform $|\mathcal{A}_+|$ compressions of $X$ yielding a non-separating $\mathbb{Z}_3$-folded surface of Euler characteristic $\chi(X)+ 2|\mathcal{A}_+| \geq \chi(X)+ |\mathcal{A}|$. $\chi_3([X])$ is maximal so it follows that $\chi_3([X]) \geq \chi(X)+ |\mathcal{A}|$, completing the proof. 
\end{proof}

It follows from applying Lemma~\ref{eq:allAdeg5} that:

\begin{equation}\label{eq:chiA}
    3\mathfrak{e} - 3 |\mathcal{A}| \geq 3\mathfrak{e} + 3\chi(X) - 3\chi_3([X])
\end{equation}

For each edge $e$, define

\begin{equation}\label{eq:dfnpsi}
\psi(e) = \begin{cases}
    s(e)  & \deg(e) = a_e  \\ 
    s(e) -3 & \text{else}
\end{cases}
\end{equation}

We call $\psi(e)$ the {\em surplus} of $e$. The total surplus of all even edges of $\mathcal{T}$ is given in Equation~\eqref{eq:psicharge}, which follows from noting that the surplus function $\psi$ applies a penalty of $3$ to $s$ of each even edge, except for those in $\mathcal{A}$. Substituting in Equation~\eqref{eq:chargeeq} yields Equation~\eqref{eq:tetrahedratypes}.

\begin{align}
    \sum_{\varphi(e)=0}\psi(e) &=\sum_{\varphi(e)=0} s(e) - 3\mathfrak{e} + 3|\mathcal{A}| \label{eq:psicharge}\\
    &= A + \frac{5}{2}B + 4C + D - 3\mathfrak{e} + 3|\mathcal{A}| \label{eq:tetrahedratypes}
\end{align}

We have the following result.

\begin{lemma}[Surplus inequality]\label{lem:psiineq}
    Consider the inequality
    \begin{equation}\label{eq:surplusinequality}
        \sum_{\varphi(e)=0}\psi(e) \geq 0
    \end{equation}
    If Inequality~\eqref{eq:surplusinequality} holds, then 
    Inequality~\eqref{eq:chargeineq} holds as well.
\end{lemma}

We call Inequality~\eqref{eq:surplusinequality} the {\em surplus inequality}. Showing that it 
holds is the content of the remainder of the proof of Theorem~\ref{thm:main}.

\begin{proof}
    Suppose that $\sum_{\varphi(e)=0}\psi(e) \geq 0$. Substituting into Equation~\eqref{eq:psicharge} yields
    \[\sum_{\varphi(e)=0} s(e) \geq 3\mathfrak{e} - 3|\mathcal{A}|\]
    Combining this inequality with Equation~\eqref{eq:chiA} produces 
    Inequality~\eqref{eq:chargeineq}:
    \[\sum_{\varphi(e)=0} s(e) \geq 3\mathfrak{e} - 3|\mathcal{A}| \geq 3\mathfrak{e} + 3\chi(X) - 3\chi_3([X])\]    
\end{proof}

\subsection{Edges of low degree}

Our goal is to prove the surplus inequality $\sum_{\varphi(e)=0}\psi(e) \geq 0$. To do this we
analyse the even edges of $\mathcal{T}$ with negative surplus, which we call the 
\textit{deficit edges}. It follows from the definition of $s$ that if $e$ is of degree at 
least $6$ then $\psi(e) \geq 0$. For each edge $e$ we can describe a \textit{model} of $e$ in 
$\mathcal{T}$ as the triangulation with the following properties:

\begin{itemize}
    \item The tetrahedra in the model of $e$ are the tetrahedra of $\mathcal{T}$ whose image in $|\mathcal{T}|$ contains $e$.
    \item The triangle pairings in the model of $e$ are the triangle pairings of $\mathcal{T}$ whose image in $|\mathcal{T}|$ contains $e$.
\end{itemize}

\textit{Models} of edges of degree less than $6$ that appear in minimal triangulations are classified in \cite[§3]{Jaco_2009}. A model of an edge $e$ of $\mathcal{T}$ is obtained by the quotient of the tetrahedra in which $e$ has a pre-image by the face pairing maps between triangles in which $e$ is an edge. We can associate to each edge a \textit{cycle type} by listing the types of the tetrahedra which contain its pre-images in the cyclic order given by the link of the edge. Lemma~\ref{lem:classifydeficitedges} classifies the deficit edges in minimal triangulations.

\begin{lemma}\label{lem:classifydeficitedges}
    The following is a complete list of all admissible cycle types and models of deficit 
    edges in a minimal, 0-efficient triangulation $\mathcal{T}$ of a closed, orientable, 
    connected, irreducible $3$-manifold 
    $M \not = L(3,1)$, with non-trivial 
    second $\mathbb{Z}_3$-homology.
    \[
        \begin{array}{c|c|c|c}
            \toprule
            \deg(e)&\text{Cycle type}&\text{Admissible models}&\psi(e)\\
            \midrule
            3&333&\mathcal{S}(1,3,4)&-1\\
            \midrule
            4&1122&X_{4;4}&-1/3\\
            4&1232&X_{4;4}&-1/6\\
            4&3333&X_{4;2}^2,X_{4;3}^0,X_{4;3}^1,X_{4;4}& -1/3\\
            \bottomrule
        \end{array}
    \]
\end{lemma}

\begin{proof}
    As mentioned before, we only need to consider edge models of degree up to $5$. In 
    \cite[§3]{Jaco_2009}, all such edge models are explicitly listed. We use their notation 
    here. Since $\mathcal{T}$ is minimal and $|\mathcal{T}|$ is not in 
    $\{S^3$, $L(3,1)$, $L(4,1)$, $L(5,1)$, $L(7,2)\}$
    (note that $M$ must have non-trivial $\mathbb{Z}_3$-homology), 
    \cite[Proposition 8]{Jaco_2009} implies that there are no edges of degree $1$ or $2$. 
    It follows from \cite[Proposition 9, p. 8]{Jaco_2009} that all edges of degree $3$ must be
    contained in a $2$-tetrahedron layered solid torus of type $\mathcal{S}(1,3,4)$. The core
    tetrahedron contains $e$ twice in opposite edges, meaning it is either of type $3$ 
    or $1$. If it is of type $3$, the second tetrahedron must also be of type $3$ and we 
    obtain the edge-model described above. If it is of type $1$, the quadrilateral is 
    identified by the face gluings to form a M\"obius band, contradicting the orientability 
    of $\text{Int}(X)$. 

    For edges of degree $4$, \cite[Proposition 10]{Jaco_2009} implies that there are four 
    distinct possible edge models $X_{4;2}^2$, $X_{4;3}^0$, $X_{4;3}^1$ and $X_{4;4}$ since 
    $|\mathcal{T}|$ is not in $\{S^3, L(4,1), L(8,3), S^3/Q_8\}$. The first three of these 
    are in layered solid tori, and the fourth is the triangulation of an octahedron with one 
    interior edge. For all of them, having only type $3$ tetrahedra is a valid solution, 
    leading to a surplus of $4 \cdot \frac{2}{3} -3 = - \frac{1}{3}$, as listed above. For 
    model $X_{4;2}^2$, we layer on top of the even degree $3$ edge in $\mathcal{S}(1,2,3)$. 
    This does not allow any other tetrahedron type than $3$, forcing the other tetrahedron 
    to be of type $3$ as well and there are no more cases to consider. For model $X_{4;3}^0$, 
    we layer two tetrahedra onto the degree $2$ edge of $\mathcal{S}(1,2,3)$. By the argument 
    from the degree $3$ case we conclude that the core tetrahedron must be of type $3$. This 
    leaves the other two tetrahedra to be either both of type $3$, or both of type $2$. The 
    former case is already listed, the latter case leads to a surplus of 
    $2 (\frac{5}{6}+ \frac{2}{3}) - 3 = 0$. For model $X_{4;3}^1$, the tetrahedron that is 
    seen twice can only be of type $1$ or $3$. If it is of type $1$, the other two tetrahedra 
    must be of type $1$ or $4$. Since $\text{Int}(X)$ is orientable, neither of them can be 
    of type $1$. If both are of type $4$, we obtain a surplus of $2 (1 + \frac{1}{2}) - 3 = 0$.
    If the tetrahedron is of type $3$, the other two tetrahedra must be of type $3$ as well, and
    we obtain cyclic type $3333$, which is already listed.
    For model $X_{4;4}$, we do not have any self-identifications. The possible cycle types are
    still constrained by the number of even pre-images in pairs of adjacent faces of tetrahedra 
    of each type. The pair of triangles sharing an even edge pre-image in a type $1$ or $4$ 
    tetrahedron always have exactly two pre-images of odd edges in them. In a type $2$ 
    tetrahedron one triangle has two pre-images of odd edges and one triangle is all-even. It 
    follows from the orientation constraints in Definition~\ref{dfn:z3foldedsurfaces} that 
    $1114$ and $1444$ are both impossible cycle types. Hence the only possible cycle types for
    a degree $4$ edge are $1111$, $1122$, $1144$, $1224$, $1232$, $1414$, $2222$, $2233$, 
    $2244$, $2324$, $3333$ and $4444$. Of these, only $1122$, $1232$ and $3333$ have negative 
    surpluses and are listed in the table above.
    
    For edges of degree $5$, a negative surplus occurs only if the cycle type is 
    $11112$, $11113$, $11133$ and $11313$. But these are all impossible configurations since 
    they would require an all-even triangle to be identified to a triangle with two odd edge 
    pre-images.    
    This concludes the case analysis.
\end{proof}

\begin{figure}[htbp]
    \centering
    \begin{tikzpicture}[
        image/.style={inner sep=0pt, outer sep=0pt},
        arrow/.style={-{Stealth}, thick, shorten <=2pt, shorten >=2pt}
    ]
        \node[image] (top1) at (0,0)
            {\includegraphics[width=0.28\linewidth]{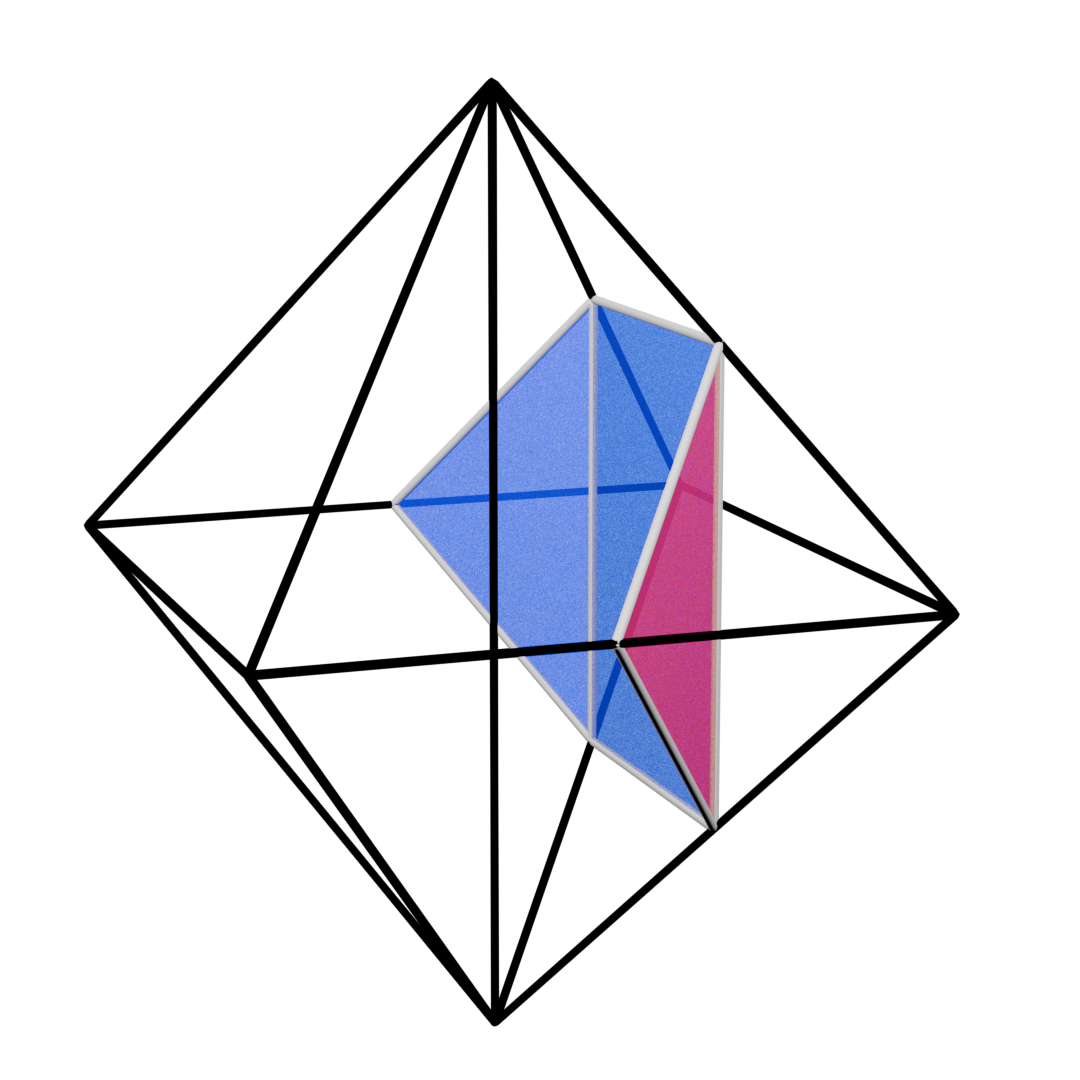}};
        \node[image] (top2) at (0.33\linewidth,0)
            {\includegraphics[width=0.28\linewidth]{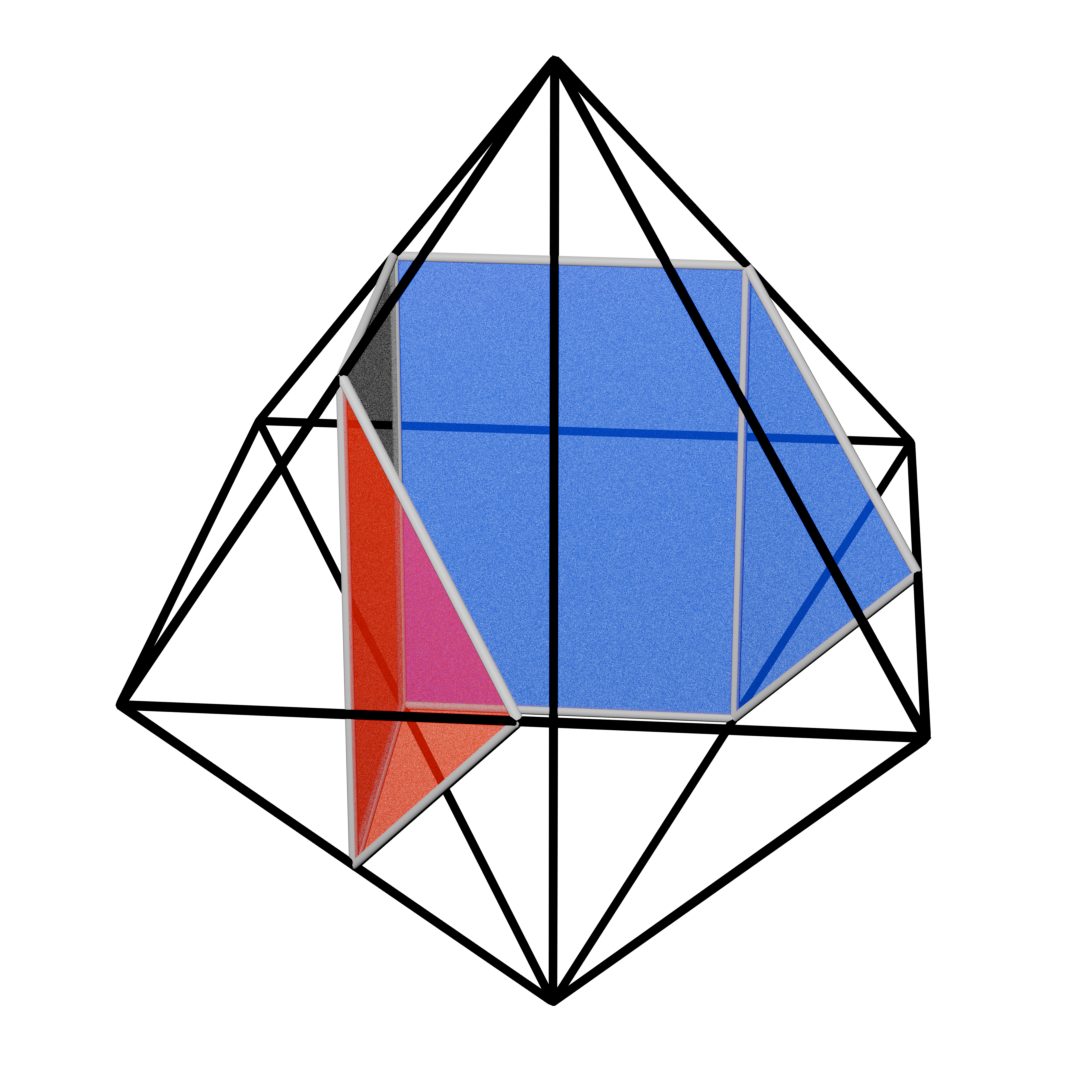}};
        \node[image] (top3) at (0.66\linewidth,0)
            {\includegraphics[width=0.28\linewidth]{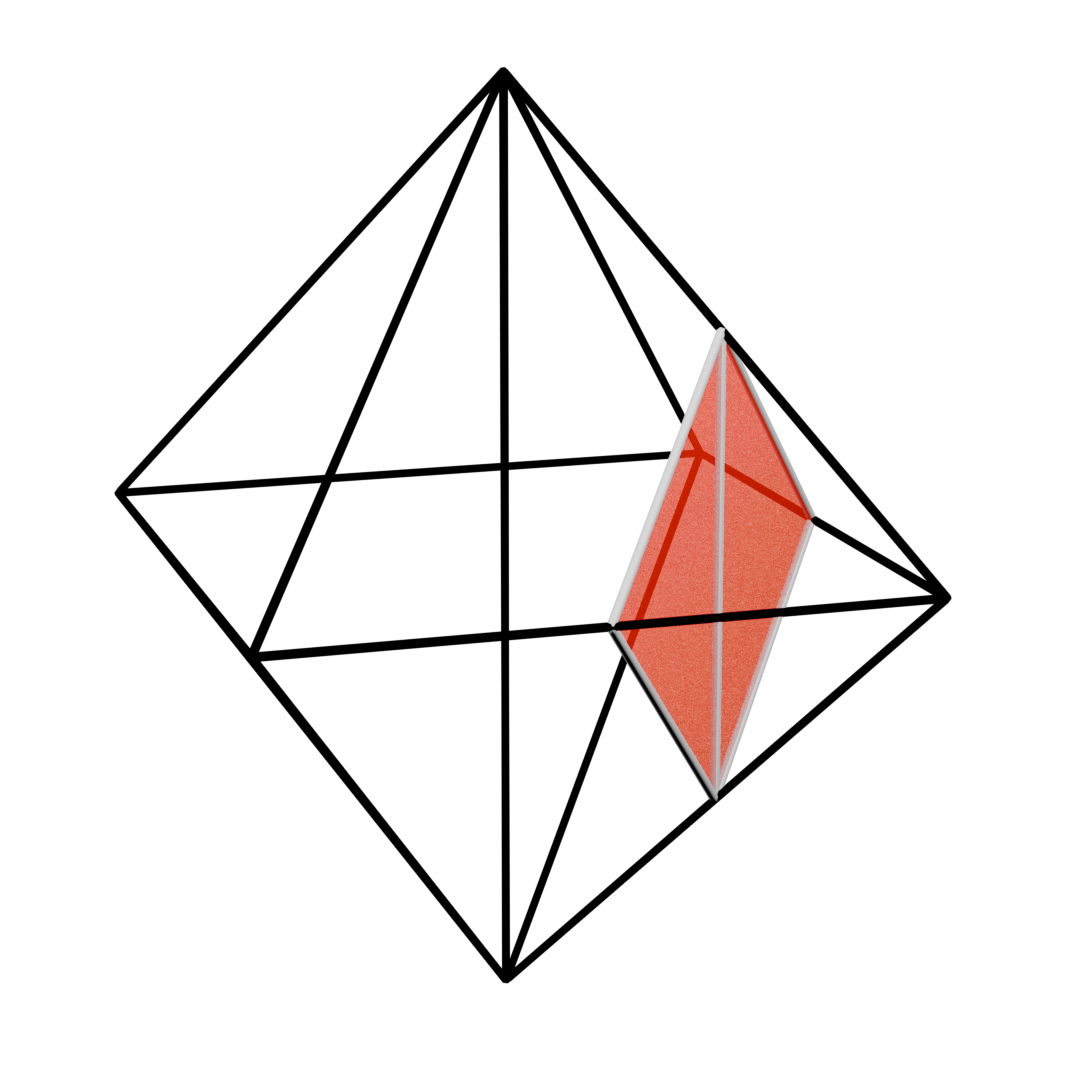}};

        \node[above=3mm of top1] {1232};
        \node[above=3mm of top2] {1122};
        \node[above=3mm of top3] {2233};

        \node[image, below=12mm of top1] (bottom1)
            {\includegraphics[width=0.28\linewidth]{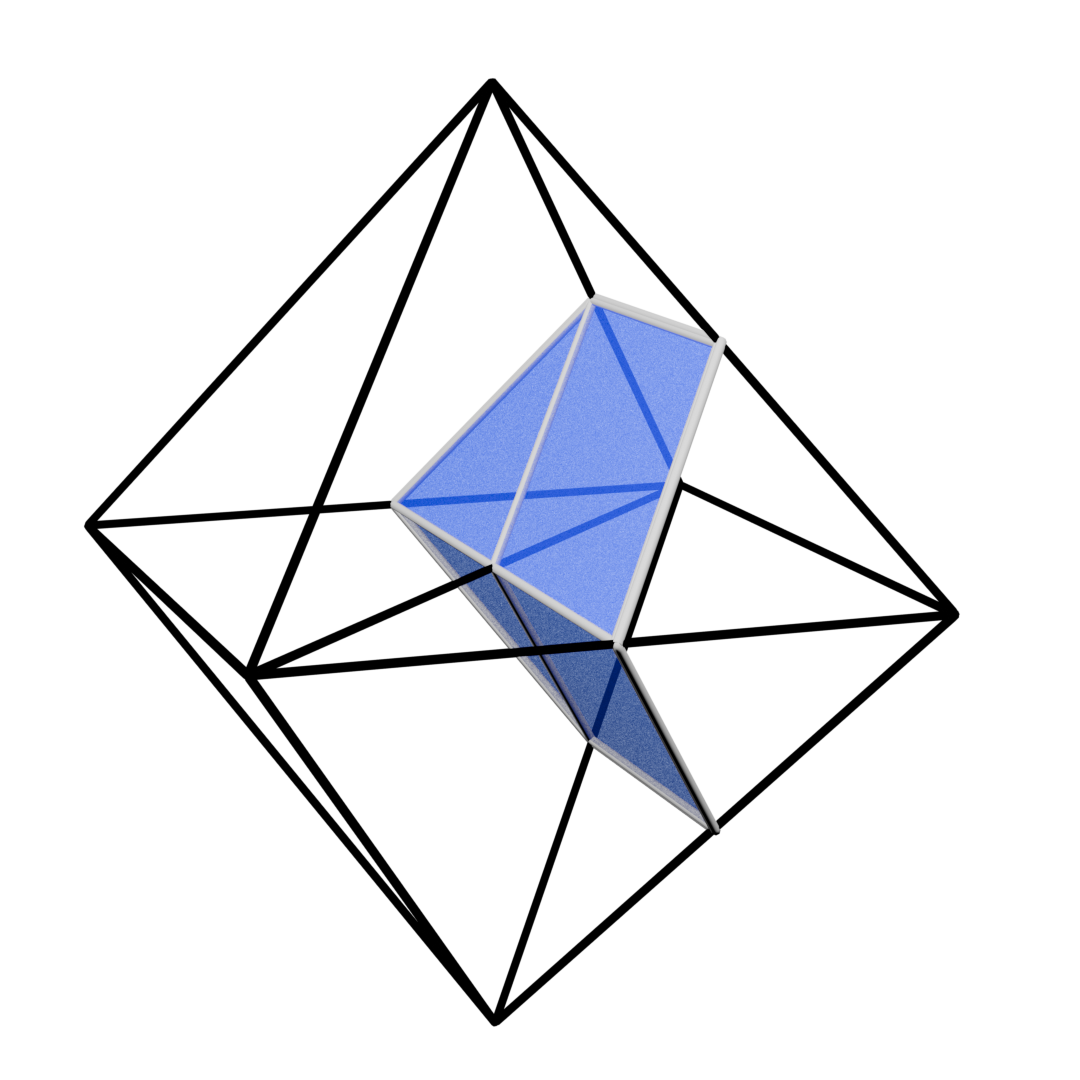}};
        \node[image, below=12mm of top2] (bottom2)
            {\includegraphics[width=0.28\linewidth]{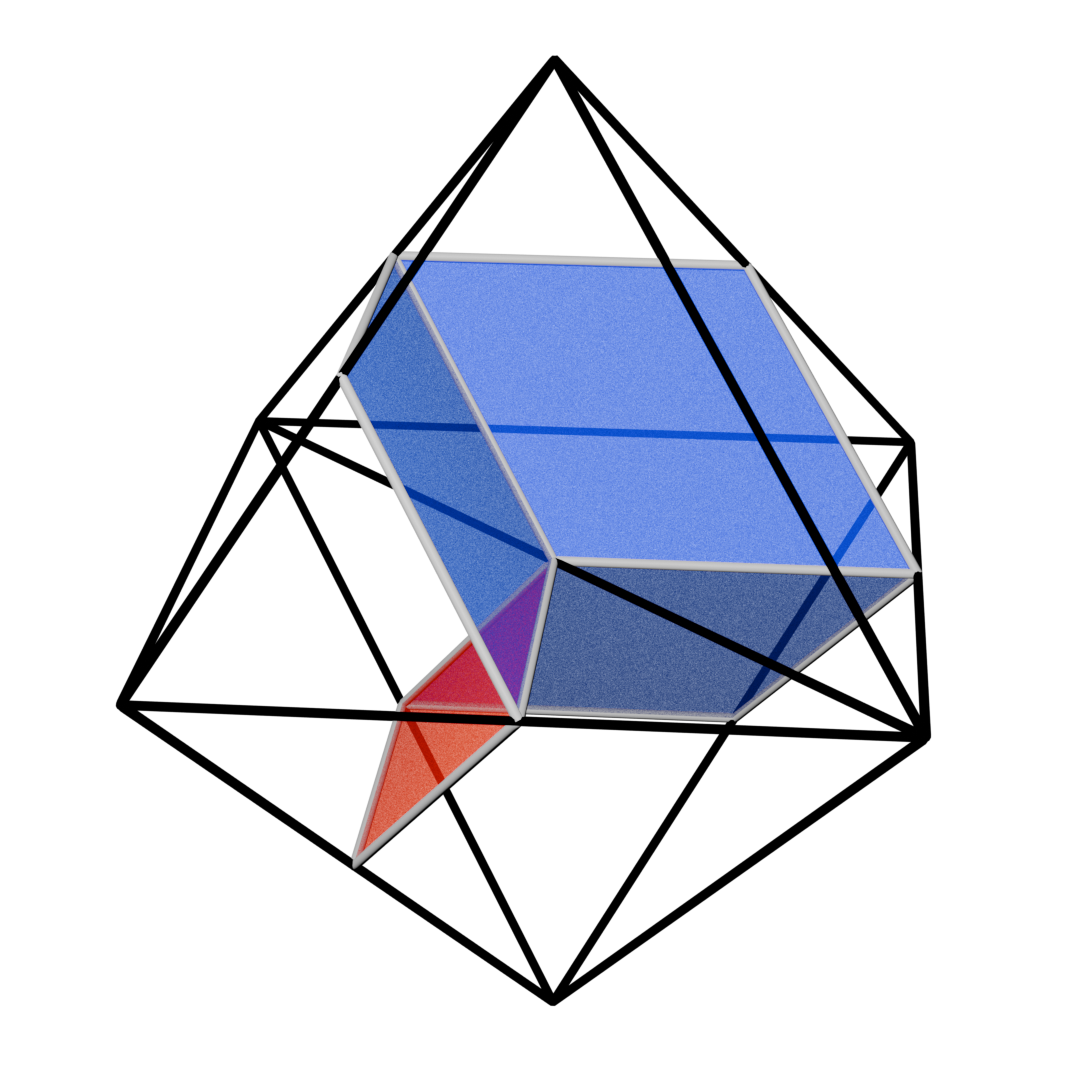}};
        \node[image, below=12mm of top3] (bottom3)
            {\includegraphics[width=0.28\linewidth]{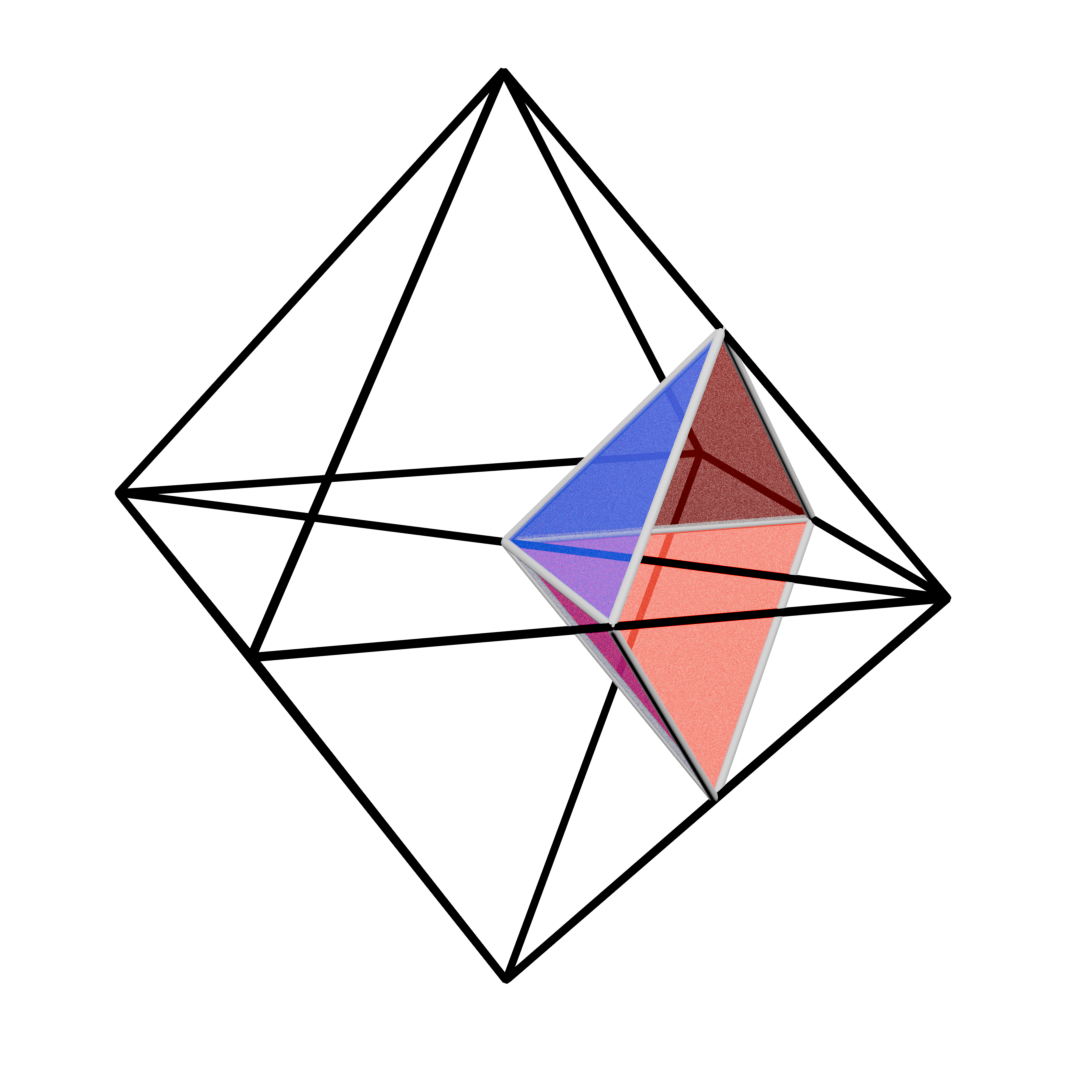}};

        \draw[arrow] (top1.south) -- (bottom1.north);
        \draw[arrow] (top2.south) -- (bottom2.north);
        \draw[arrow] (top3.south) -- (bottom3.north);
    \end{tikzpicture}
    \caption{$(4-4)$-Pachner moves applied to tetrahedra with a central even edge of cycle types $1122$, $1232$ and $2233$ that result in another minimal triangulation with strictly fewer even edges since the new diagonal edge is now odd.}
    \label{fig:44movereplacements}
\end{figure}

Cycles of type $1232$ and $1122$ occur only in the four-tetrahedron triangulation of the
octahedron model $X_{4;4}$. A $(4-4)$-move replacing the degree $4$ edge with a different
diagonal of the octahedron that intersects the canonical $\mathbb{Z}_3$-folded surface yields 
another triangulation with the same number of tetrahedra and strictly fewer even edges as in 
Figure~\ref{fig:44movereplacements}.
Hence after a sufficient number of replacements we can assume that there are no even edges 
with these cycle types in the triangulation. This same procedure can also be applied to 
eliminate any even degree $4$ edges of type $2233$ in an octahedron configuration. Even 
though the even edge in this model contributes a surplus of zero, this simplifies the 
proof in Section~\ref{ssec:subcomplex}.
Hence, we obtain a triangulation in which every even edge $e$ that contributes a deficit 
($\psi(e)<0$) has the property that $c_e = \deg(e) \in \{3,4\}$. 
For the remainder of this article, we assume that $\mathcal{T}$ satisfies this condition.

\subsection{Subcomplex $\mathcal{C}$ and triangulation $\tilde{\mathcal{C}}$}
\label{ssec:subcomplex}

We define the \emph{type $3$ subcomplex}, denoted $\mathcal{C}$, to be the subcomplex of $\mathcal{T}$ spanned by the type $3$ tetrahedra in the triangulation. Every edge in this subcomplex is even, and the deficit edges of $\mathcal{T}$ are all contained in the interior of $\mathcal{C}$. To see this note that only edges in the interior of $\mathcal{C}$ satisfy $c_e = \deg (e)$. We seek to show that the total surplus of all the even edges of $\mathcal{C}$ is non-negative. We define the \emph{type $3$ subtriangulation}, denoted $\tilde{\mathcal{C}}$, to be the triangulation whose tetrahedra are the tetrahedra of $\mathcal{T}$ contained in $\mathcal{C}$ and whose triangle pairings are the triangle pairings in $\mathcal{T}$ between faces of tetrahedra in $\mathcal{C}$. This construction induces a map $|\tilde{\mathcal{C}}| \rightarrow |\mathcal{C}|$ given by the quotient of the map that sends the tetrahedra of $\tilde{\mathcal{C}}$ to the corresponding tetrahedra of $\mathcal{T}$. This map restricts to an embedding on the interior of $|\tilde{\mathcal{C}}|$.

Define the following:
\begin{itemize}
    \item Let $k$ be the number of tetrahedra in $\tilde{\mathcal{C}}$.
    \item Let $e_{\partial}$ and $f_{\partial}$ be the number of edges and triangles in $\partial |\tilde{\mathcal{C}}|$, respectively.
    \item Let $e_{\text{int}}$ and $f_{\text{int}}$ be the number of edges and triangles in $\text{Int}(|\tilde{\mathcal{C}}|)$, respectively.
    \item Let $v$ be the number of vertices on the boundary of $\tilde{\mathcal{C}}$.
    \item Let $d_1$ be the number of edges in $\partial |\tilde{\mathcal{C}}|$ of degree $1$.
    \item Let $d_2$ be the number of edges in $\partial |\tilde{\mathcal{C}}|$ of degree $2$ that are not on the boundary of a component that is the layered solid torus $\mathcal{S}(1,2,3)$.
\end{itemize}

Lemma~\ref{lem:sufficientconditionpinchedhandle} gives a sufficient condition expressed in terms of the combinatorics of $\tilde{\mathcal{C}}$ for the global surplus $\sum_{\varphi(e)=0}\psi(e)$ to be non-negative.

\begin{lemma}\label{lem:sufficientconditionpinchedhandle}
    If 
    \begin{equation}\label{eq:pinchedhandleineq}
        4k - 3(e_{\partial} + e_{\text{int}}) + \frac{5}{2} f_{\partial} + d_1 + d_2/2 \geq 0
    \end{equation}
    holds for $\tilde{\mathcal{C}}$, then the surplus inequality 
    $\sum_{\varphi(e)=0}\psi(e) \geq 0$ holds as well.    
\end{lemma}

\begin{proof}
    Let $E_{3}$ be the set of even edges $e$ in $\mathcal{T}$ with the property that 
    $c_e\geq 1$. After applying the $(4-4)$-Pachner moves in 
    Figure~\ref{fig:44movereplacements}, it follows from 
    Lemma~\ref{lem:classifydeficitedges} that every deficit edge is in $E_3$ and we can 
    obtain a lower bound on the global surplus $\sum_{\varphi(e)=0}\psi(e)$ by 
    restricting the sum to $E_3$ in
    \begin{equation}\label{eq:E3lowerbound}
        \sum_{\varphi(e)=0}\psi(e) \geq \sum_{e\in E_3}\psi(e) = \sum_{e\in E_3}s(e) - 3|E_3|
    \end{equation}
    Every edge in $E_3$ is in the image of the map $|\tilde{\mathcal{C}}|\rightarrow |\mathcal{T}|$, but some of the edges in $\partial |\tilde{\mathcal{C}}|$ may be identified together by triangle pairings of $\mathcal{T}$ not present among the pairings of $\tilde{\mathcal{C}}$, leading to
    \begin{equation}\label{eq:E3upperbound}
        e_{\text{int}}+e_{\partial} \geq |E_{3}|.
    \end{equation}
    All six edges on the boundary of each type $3$ tetrahedron are pre-images of edges in $E_3$ and every type $3$ tetrahedron has a unique pre-image in $\tilde{\mathcal{C}}$ so we have 
    \begin{equation}\label{eq:countingtype3}
        \sum_{e\in E_3}c_e = 6C = 6k.
    \end{equation}
    Every boundary triangle in $\partial(|\tilde{\mathcal{C}}|)$ is mapped to the boundary of a unique type $2$ tetrahedron of $\mathcal{T}$ and all three even edges of that type $2$ tetrahedron are pre-images of edges in $E_3$, yielding
    \begin{equation}\label{eq:countingtype2}
        \sum_{e\in E_3}b_e \geq 3f_{\partial}.
    \end{equation}
    Applying Inequalities~\eqref{eq:E3lowerbound}, \eqref{eq:E3upperbound}, 
    \eqref{eq:countingtype2} and Equation~\eqref{eq:countingtype3} produces a na\"ive 
    lower bound on the global surplus in terms of the combinatorics of 
    $\tilde{\mathcal{C}}$.
    \begin{align}
        \sum_{\varphi(e)=0}\psi(e) &\geq \sum_{e\in E_3}s(e) - 3|E_3|\notag\\
        &\geq \sum_{e\in E_3}s(e) - 3(e_{\text{int}}+ e_{\partial})\notag\\
        &\geq \sum_{e\in E_3}\left[\frac{5b_e}{6} + \frac{2c_e}{3}\right] - 3(e_{\text{int}}+ e_{\partial}) \notag\\
        &\geq \frac{5}{6}\cdot 3f_{\partial} + \frac{2}{3}\cdot 6k - 3(e_{\text{int}}+ e_{\partial}) \notag\\
        &= \frac{5}{2}f_{\partial}+4k- 3(e_{\text{int}}+ e_{\partial})\label{eq:naivelowerbound}
    \end{align}
    For each edge $e$ in $E_3$ let $m_e$ be the number of pre-images of $e$ in $\tilde{\mathcal{C}}$. The gap in Inequality~\eqref{eq:E3upperbound} is exactly $\sum_{e\in E_3}\left[m_e-1\right]$, introducing a gap of $3\sum_{e\in E_3}\left[m_e-1\right]$ in Inequality~\eqref{eq:naivelowerbound}. Now suppose we aggregate over the edges of $\tilde{\mathcal{C}}$ instead. If we share this gap equally among the edges in $\tilde{\mathcal{C}}$ identified with $e$ we have a gap of $\frac{3(m_e-1)}{m_e}$ for each edge above $e$. If $m_e>1$ then $\frac{3(m_e-1)}{m_e} \geq \frac{3}{2}$. In particular, for any edge in $\tilde{\mathcal{C}}$ with the property that there is at least one other edge in $\tilde{\mathcal{C}}$ that shares the same image in $\mathcal{T}$ we can add $\frac{3}{2}$ to the right hand side of Inequality~\eqref{eq:naivelowerbound} without breaking it. 

    Let $e$ be an edge in $\partial|\tilde{\mathcal{C}}|$ that is of degree $1$ in $\tilde{\mathcal{C}}$. If there is any other edge in $\partial|\tilde{\mathcal{C}}|$ that shares the same image in $E_3$ we attribute a gap of at least $\frac{3}{2}$ to it as per the preceding argument. Now assume that $e$ is the unique pre-image of an edge in $E_3$. Inequality~\eqref{eq:naivelowerbound} estimates a lower bound on $s(e)$ by counting the contribution coming from the single pre-image of $e$ in a type $3$ tetrahedron and the two pre-images of $e$ in the type $2$ tetrahedra identified to the pair of triangles of the type $3$ tetrahedron that meet at $e$. We can sharpen Inequality~\eqref{eq:naivelowerbound} by adding the contribution from any remaining pre-images of $e$ in $\mathcal{T}$. Since $232$ is not an admissible configuration of tetrahedra in a minimal triangulation, $e$ must be of degree at least $4$ in $\mathcal{T}$. We have already eliminated models $1232$, so if $e$ has degree four the only possibility is $2324$ and the missing contribution from the pre-image in the type $4$ tetrahedron is $1$. If the degree of $e$ is greater than four then there are at least two tetrahedra unaccounted for and so the missing contribution is at least $2 \cdot \frac12 = 1$ in this case as well. In all cases the gap associated to an edge of degree $1$ in $\partial|\tilde{\mathcal{C}}|$ is at least $1$. By definition, there are $d_1$ such edges.
    
    Similarly, let $e$ be an edge in $\partial|\tilde{\mathcal{C}}|$ that is of degree $2$ in $\tilde{\mathcal{C}}$ which is not on the boundary of a single tetrahedron layered solid torus $\mathcal{S}(1,2,3)$. If there is any other edge in $\partial|\tilde{\mathcal{C}}|$ that shares the same image in $E_3$ we again attribute a gap of at least $\frac{3}{2}$ to it. Now assume that $e$ is the unique pre-image of an edge in $E_3$. It follows that the cycle type of the image of $e$ in $\mathcal{T}$ has a single block of two type $3$ tetrahedra in it with a type $2$ tetrahedron on each side. If the degree of the image of $e$ in $\mathcal{T}$ is four then we have a local model of type $2233$. We have already excluded the octahedral model with a $2233$ cycle type via a $(4-4)$-Pachner move, so this model must be $X_{4;3}^0$ and the connected component of $\tilde{\mathcal{C}}$ containing $e$ is the single tetrahedron layered solid torus $\mathcal{S}(1,2,3)$ which contradicts our initial assumption. Hence the degree of the image of $e$ in $\mathcal{T}$ is at least $5$ and there exists at least one missing tetrahedron with a contribution of at least $1/2$. In all cases the gap associated to an edge of degree $2$ in $\partial|\tilde{\mathcal{C}}|$ not on the boundary of an $\mathcal{S}(1,2,3)$ component is at least $\frac{1}{2}$. By definition there are $d_2$ such edges.

    Together these observations imply that there is a gap of at least $d_1+d_2/2$ in
    Inequality~\eqref{eq:naivelowerbound} and so:

    \[\sum_{\varphi(e)=0}\psi(e) \geq 4k - 3(e_{\partial} + e_{\text{int}}) + \frac{5}{2} f_{\partial} + d_1 + d_2/2\]

    The statement of the lemma follows immediately.
    
\end{proof}

Lemma~\ref{lem:Ccomplextopology} describes the topology of $|\tilde{\mathcal{C}}|$. 

\begin{lemma}\label{lem:Ccomplextopology}
    Let $\tilde{\mathcal{C}}$ be the type $3$ sub-triangulation derived from a minimal 
    triangulation $\mathcal{T}$ of a lens space $L(3n,1)$, $n\geq 2$, and associated 
    canonical $\mathbb{Z}_3$-folded surface $X$. Each component $H$ of 
    $\tilde{\mathcal{C}}$ has the following properties:
    
    \begin{itemize}
        \item $\text{Int}(H)$ is an orientable $3$-manifold.
        \item Every vertex of $H$ lies on its boundary and the link of every vertex is a connected, compact, planar $2$-manifold with boundary.
        \item Every interior edge of $H$ has degree at least $3$.
        \item If we let $\mathring{H}$ denote the truncation of $H$ at all its vertices, $\partial \mathring{H}$ must be a connected, closed, orientable $2$-manifold.
    \end{itemize}
    
\end{lemma}

\begin{proof}
    $\tilde{\mathcal{C}}$ is a triangulation with triangle pairings that are a subset of the triangle pairings in $\mathcal{T}$ and $|\mathcal{T}|$ is a closed $3$-manifold. It follows that $\tilde{\mathcal{C}}$ is a valid triangulation and so every point has the local topology of a $3$-manifold, except possibly the vertices. Since $\mathcal{T}$ has a single vertex, every vertex of $\tilde{\mathcal{C}}$ must lie on its boundary, which is non-trivial since $\varphi \neq 0$ and $X\neq \emptyset$. It follows that $\text{Int}(|\tilde{\mathcal{C}}|)$ is an orientable $3$-manifold. The link of each vertex must be connected since the corners of tetrahedra in $\tilde{\mathcal{C}}$ represent the same vertex only when there exists a sequence of triangle pairings relating them. Moreover, since these triangle pairings are a subset of the triangle pairings in $\mathcal{T}$ and $|\mathcal{T}|$ is a closed $3$-manifold, it follows that the link of each vertex of $\tilde{\mathcal{C}}$ embeds in $S^2$. Since $\mathcal{T}$ is a minimal triangulation of a manifold other than $S^3, L(3,1)$ and $L(4,1)$, it follows from \cite[Proposition 8]{Jaco_2009} that every edge must have degree at least $3$. Edges in the interior of $|\tilde{\mathcal{C}}|$ have the same degree as their image in $|\mathcal{T}|$ so they also must have degree at least $3$.

    Let $H$ be a connected component of $|\tilde{\mathcal{C}}|$. Since $H$ has the topology of a compact, orientable $3$-manifold everywhere except for at its vertices, it follows that truncating these vertices to form $\mathring{H}$ yields a connected, compact, orientable $3$-manifold. The only property left to check is that $\partial \mathring{H}$ is connected. 
    
    We assume for a contradiction that $\partial \mathring{H}$ has more than one component and 
    let $\Sigma$ be one of them. Push $\Sigma$ slightly into the interior of $H$ to obtain a 
    parallel closed surface separating $\Sigma$ from the remaining boundary components. 
    The map $|\tilde{\mathcal{C}}|\rightarrow |\mathcal{T}|$ restricts to an embedding on the
    interior of $H$ and so it carries the parallel surface to a corresponding closed orientable
    surface in $|\mathcal{T}|$. This map also carries each boundary triangle of $H$ to a triangle in $|\mathcal{T}|$ on the boundary of a unique type $2$ tetrahedron. We use this to locate $X$ on both sides of the surface. Since $\mathcal{T}$ is a triangulation of a lens space, every embedded closed orientable surface separates the ambient space so we deduce that $X$ has multiple components. This contradicts \cite[Lemma 4.8, p. 28]{hobkirk2026} which implies that $X$ is connected since lens spaces have second $\mathbb{Z}_3$-homology of rank at most $1$. Hence $\partial \mathring{H}$ is connected. 

\end{proof}

\subsection{Component-wise verification}

For each connected component $H$ of $\tilde{\mathcal{C}}$ let $k^H$ be the number of 
tetrahedra in $H$. Define $e_{\partial}^H$, $e_{\text{int}}^H$, $f_{\partial}^H$, $v^H$, 
$d_1^H$ and $d_2^H$ similarly. 
The quantities $k,e_{\text{int}},e_{\partial},f_{\partial},d_1,d_2$, on which 
Inequality~\eqref{eq:pinchedhandleineq} depends, can be decomposed into a sum of the 
contributions $k^H,e_{\text{int}}^H,e_{\partial}^H,f_{\partial}^H,d_1^H,d_2^H$ from 
each connected component $H$ of $\tilde{\mathcal{C}}$. We consider the following inequalities:
\begin{equation}\label{eq:s123case}
        4k^H - 3(e_{\partial}^H + e_{\text{int}}^H) + \frac{5}{2} f_{\partial}^H + d_1^H \geq 0
\end{equation}
\begin{equation}\label{eq:generalcase}
    4k^H - 3(e_{\partial}^H + e_{\text{int}}^H) + \frac{5}{2} f_{\partial}^H + d_1^H + d_2^H/2 \geq 0
\end{equation}
It follows that it is sufficient 
to prove that Inequality~\eqref{eq:s123case} holds when $H$ is the single tetrahedron 
layered solid torus $\mathcal{S}(1,2,3)$ and Inequality~\eqref{eq:generalcase} holds 
in every other case. The first claim follows immediately from the observation that 
$\mathcal{S}(1,2,3)$ has one tetrahedron, three boundary edges of degrees $1,2$ and 
$3$, no interior edges and two boundary triangles and hence the LHS of 
Inequality~\eqref{eq:s123case} evaluates to $1$. 
The second claim is proved in Proposition~\ref{prop:desiredbound}.

\begin{proposition}\label{prop:desiredbound}
    Inequality~\eqref{eq:generalcase} holds for every component $H$ of~$\tilde{\mathcal{C}}$.
\end{proposition}

\begin{proof}

    We start with Lemma~\ref{lem:Ccomplextopology}. Let $\mathring{H}$ be the orientable $3$-manifold with connected boundary resulting from truncation of all the vertices of $H$. Let $g^H$ be the genus of $\partial \mathring{H}$. For each vertex $x$ of $H$ let $b_x$ be the number of boundary components in the link of $x$. 
        
    Counting the number of edges on the boundary of each triangle and the edges identified by gluings results in
    
    \begin{equation}
        \frac{3}{2}f_{\partial}^H = e_{\partial}^H\label{eq:bdryfaces}
    \end{equation}
    
    A similar argument counting faces on the boundary of each tetrahedron gives $4k^H = 2f_{\text{int}}^H + f_{\partial}^H.$ This rearranges to $f_{\text{int}}^H = 2k^H - f_{\partial}^H/2$, and substituting Equation~\eqref{eq:bdryfaces} for $f_{\partial}^H$ gives
    
    \begin{equation}
        f_{\text{int}}^H = 2k^H - e_{\partial}^H/3 \label{eq:intfaces}
    \end{equation}

    $H$ is obtained from its truncation $\mathring{H}$ by coning off each vertex link. The cone on each vertex link has Euler characteristic $1$ and the link of a vertex $x$ has Euler characteristic $2- b_x$. It follows that
    
    \begin{align}
        \chi(H) &= \chi(\mathring{H})+ \sum_{x \in H^{(0)}} \left[1 - (2-b_x)\right]\notag\\
        &= \chi(\mathring{H}) + \sum_{x \in H^{(0)}} \left[b_x-1\right]\label{eq:boundary}\\
        &= 1-g^H + \sum_{x \in H^{(0)}} \left[b_x-1\right] \label{eq:genusvalue}
    \end{align}

    Now compare $\partial \mathring{H}$ to $\partial H$. The former is obtained by filling the boundary circles in the link of each vertex $x$ with a connected planar $2$-manifold of Euler characteristic $2-b_x$. The latter is obtained by filling the same boundary circles with their cone to the vertex $x$ which has Euler characteristic $1$. It follows that

    \[\chi(\partial H) = \chi(\partial \mathring{H}) + \sum_{x \in H^{(0)}} \left[b_x-1\right]\]

    Since $\mathring{H}$ is an orientable $3$-manifold with boundary it follows that $\chi(\partial \mathring{H}) = 2\chi(\mathring{H})$ and so Equation~\eqref{eq:boundary} can be rearranged into

    \begin{equation}\label{eq:boundary2}
        2\chi(H) = \chi(\partial H) + \sum_{x \in H^{(0)}} \left[b_x-1\right].
    \end{equation}

    Cell counts for $H$ and $\partial H$ yield the following, noting that all vertices of $H$ lie on the boundary

    \begin{align}
        \chi(H) &= -k^H + f_{\partial}^H + f_{\text{int}}^H - e_{\partial}^H - e_{\text{int}}^H + v^H \label{eq:Hcellcount}\\
        \chi(\partial H) &= f_{\partial}^H- e_{\partial}^H + v^H. \notag
    \end{align}

    Substituting these counts into Equation~\eqref{eq:boundary2} yields

    \[-2k^H + f_{\partial}^H + 2f_{\text{int}}^H - e_{\partial}^H - 2e_{\text{int}}^H + v^H = \sum_{x \in H^{(0)}} \left[b_x-1\right].\]

    Using $4k^H = 2f_{\text{int}}^H + f_{\partial}^H$ we obtain 

    \begin{equation}\label{eq:verts}
        v^H = -2k^H + e_{\partial}^H + 2e_{\text{int}}^H + \sum_{x \in H^{(0)}} \left[b_x-1\right].
    \end{equation}

    Substituting this back into Equation~\eqref{eq:Hcellcount} yields \eqref{eq:subback}; substituting Equations \eqref{eq:bdryfaces} and \eqref{eq:intfaces} yields \eqref{eq:subbdryint}; 

    \begin{align}
        \chi(H) &= -3k^H + f_{\partial}^H + f_{\text{int}}^H + e_{\text{int}}^H + \sum_{x \in H^{(0)}} \left[b_x-1\right] \label{eq:subback} \\
        &= -k^H + \frac{1}{3}e_{\partial}^H  + e_{\text{int}}^H + \sum_{x \in H^{(0)}} \left[b_x-1\right] \label{eq:subbdryint}.
    \end{align}

    Multiplying \eqref{eq:subbdryint} by $4$ and rearranging to obtain the LHS of Equation~\eqref{eq:generalcase} yields \eqref{eq:lhsforce}; applying Equation~\eqref{eq:bdryfaces} this simplifies to \eqref{eq:simplifylhs}; substituting in \eqref{eq:genusvalue} this simplifies further to \eqref{eq:fullysimplified}.

    \begin{align}
        &4k^H -3(e_{\partial}^H + e_{\text{int}}^H) + \frac{5}{2}f_{\partial}^H + d_1^{H} + d_2^{H}/2 \notag\\
        &{}= -4\chi(H)+ \frac{5}{2}f_{\partial}^H - \frac{5}{3}e_{\partial}^H + e_{\text{int}}^H + 4\sum_{x \in H^{(0)}} \left[b_x-1\right] + d_1^{H} + d_2^{H}/2\label{eq:lhsforce}\\
        &{}= -4\chi(H)+ e_{\text{int}}^H + 4\sum_{x \in H^{(0)}} \left[b_x-1\right] + d_1^{H} + d_2^{H}/2\label{eq:simplifylhs} \\
        &{}= 4g^H - 4 + e_{\text{int}}^H + d_1^{H} + d_2^{H}/2 \label{eq:fullysimplified}
    \end{align}

    Inequality~\eqref{eq:generalcase} holds if and only if 

    \begin{equation}\label{eq:equivalentcondition}
        4g^H - 4 + e_{\text{int}}^H + d_1^{H} + d_2^{H}/2 \geq 0.
    \end{equation}
    
    This is trivially true whenever $g^H>0$ or $e_{\text{int}}^H\geq 4$ and so for the remainder
    of the proof we assume that $\partial \mathring{H}$ is a $2$-sphere and $H$ has fewer than
    $4$ interior edges.

    It follows from counting edge pre-images of the tetrahedra in $H$ that
    
    \begin{align}
    6k^H  &\geq  \sum_{i=1}^{\infty} i \cdot d_i^H + 3e_{\text{int}}^H \nonumber \\
    &\geq 3\,(e_{\partial}^H - d_1^{H} - d_2^{H}) + d_1^{H}+ 2d_2^{H} + 3e_{\text{int}}^H \nonumber \\
    &= 3e_{\partial}^H + 3e_{\text{int}}^H -2d_1^{H} - d_2^{H} \label{eq:tet_edge_bound}
    \end{align}

    Now recall Equation~\eqref{eq:verts} which implies that 
    \[6k^H = 3e_{\partial}^H +6e_{\text{int}}^H - 3v^H + 3\sum_{x \in H^{(0)}} \left[b_x-1\right].\] 
    
    Substituting this into Inequality~\eqref{eq:tet_edge_bound} gives
    \[3e_{\partial}^H + 6e_{\text{int}}^H - 3v^H + 3\sum_{x \in H^{(0)}} \left[b_x-1\right]\geq 3e_{\partial}^H + 3e_{\text{int}}^H -2d_1^{H} - d_2^{H}\]
    Rearranging and dividing by $3$ yields Inequality~\eqref{eq:vbound}.

    \begin{equation}\label{eq:vbound}
        e_{\text{int}}^H + d_1^{H} + \frac{d_2^{H}}{2} \geq v^H - \sum_{x \in H^{(0)}} \left[b_x-1\right] + \frac{d_1^{H} + \frac{d_2^{H}}{2}}{3}
    \end{equation}

    Now let $N$ be the \emph{normalisation} of $\partial H$, obtained by resolving all its singular vertices. 
    $N$ is homeomorphic to the surface obtained from $\partial \mathring{H}$ by replacing the link of each vertex $x$ with $b_x$ discs capping off each boundary circle in the link. Since $\partial\mathring{H}$ is a $2$-sphere, discarding the planar link of any vertex produces a planar $2$-manifold with boundary and identifying a disc to a boundary component of a planar $2$-manifold with boundary yields another planar $2$-manifold with boundary or a $2$-sphere. It follows that $N$ must be a disjoint union of $2$-spheres. The Euler characteristic of the link of a vertex $x$ is $2-b_x$ and each capping disc increases the Euler characteristic by $1$ so

    \begin{align*}
        \chi(N) &= \chi(\partial\mathring{H}) +\sum_{x \in H^{(0)}} \left[b_x - (2-b_x)\right] \\
        &= 2 + 2\sum_{x \in H^{(0)}} \left[b_x -1\right]
    \end{align*}

    Hence $N$ is a disjoint union of $1 + \sum_{x \in H^{(0)}} \left[b_x -1\right]$ spheres. Any triangulation of a $2$-sphere has at least three vertices so the induced triangulation of $N$ has at least $3+ 3\sum_{x \in H^{(0)}} \left[b_x -1\right]$ vertices. Each vertex $x$ of $\partial H$ resolves into $b_x$ vertices in the induced triangulation of $N$. It follows that $N$ has exactly $v^H + \sum_{x \in H^{(0)}} \left[b_x -1\right]$ vertices. Combining these observations we obtain the following inequality

    \[v^H - \sum_{x \in H^{(0)}} \left[b_x -1\right] \geq 3 + \sum_{x \in H^{(0)}} \left[b_x -1\right].\]

    Substituting this into the RHS of Inequality~\eqref{eq:vbound} yields

    \[e_{\text{int}}^H + d_1^{H} + \frac{d_2^{H}}{2} \geq 3 + \sum_{x \in H^{(0)}} \left[b_x -1\right] + \frac{d_1^{H} + \frac{d_2^{H}}{2}}{3}\]

    In particular if $\sum_{x \in H^{(0)}} \left[b_x -1\right] >0$ it follows immediately that $e_{\text{int}}^H + d_1^{H} + \frac{d_2^{H}}{2} \geq 4$ and Inequality~\eqref{eq:equivalentcondition} holds.

    We now assume the remaining case where $\sum_{x \in H^{(0)}} \left[b_x -1\right]= 0$. It follows that the link of every vertex of $H$ is a disc and $\partial H$ is a $2$-sphere. If $v^H\geq 4$, Inequality~\eqref{eq:vbound} implies that $e_{\text{int}}^H + d_1^{H} + \frac{d_2^{H}}{2} \geq 4$ and so Inequality~\eqref{eq:equivalentcondition} holds. We therefore also assume that $v^H<4$. 
    
    Since $\partial H$ is a triangulation of $S^2$ it follows that $v^H - e_{\partial}^H + f_{\partial}^H = 2$. Moreover substituting in Equation~\eqref{eq:bdryfaces} gives $v^H = 2 + f_{\partial}^H/2$. $\partial H$ must have at least two boundary faces so $v^H\geq 3$. It follows that $\partial H$ is a triangulation of the $2$-sphere with exactly three vertices, three edges and two triangles. It also follows from the irreducibility of lens spaces that $H$ must be a $3$-ball, but the proof does not require this. All that remains is to show that for all such triangulations $e_{\text{int}}^H + d_1^{H} + \frac{d_2^{H}}{2} \geq 4$.
    
    We now give the reader two options for how to finish the proof, one computational and one combinatorial.

    \subsubsection*{The computational option}
    
    The computational method is to exhaustively enumerate the finite family of orientable $3$-manifold triangulations extending triangulations of the $2$-sphere with two triangles. The restriction to a finite census follows from the fact that any triangulation $H$ extending the $2$-sphere with $e_{\text{int}}^H<4, v^H=3, e_{\partial}^H=3,f_{\partial}^H = 2$ and no interior vertices must have fewer than four tetrahedra. To see this note that every triangle of $H$ except the two boundary triangles is contained in exactly $2$ tetrahedra, and so $f_{\text{int}}^H+f_{\partial}^H = 2k^H+1$. Since $\chi(H) = \chi(\partial H)/2 = 1$ the Euler characteristic yields

    \[v^H - (e_{\text{int}}^H + e_{\partial}^H) + (f_{\text{int}}^H + f_{\partial}^H) - k^H   = 1.\]

    In particular, we have

    \[3 - (e_{\text{int}}^H+3) + 2k^H + 1 - k^H   = 1.\]

    This simplifies to $k^H = e_{\text{int}}^H <4$. This means that it is sufficient to enumerate all orientable $3$-manifold triangulations with at most $3$ tetrahedra, exclude any triangulations that have either an interior edge of degree less than $3$ or an interior vertex, and then verify that all remaining triangulations satisfy Inequality~\eqref{eq:equivalentcondition}. This is a relatively straightforward computation that we leave to the interested reader.

    \subsubsection*{The combinatorial option}

    For a combinatorial proof, we show directly that every triangulation of a $3$-manifold extending a triangulation of a $2$-sphere with two boundary triangles, no interior vertices and fewer than four interior edges, each of degree at least $3$, satisfies $d_1^{H} + d_2^{H}/2 + e_{\text{int}}^H \geq 4$. 
    
    There are only two combinatorially distinct triangulations of the $2$-sphere with two triangles, planar projections of which are reproduced in Figure~\ref{fig:2spheretriangulations}. We label these triangulations in terms of the \emph{boundary degree} of the three vertices, which we define to be the number of endpoints of the edges in $\partial H$ identified with each vertex. In the $(2,2,2)$ triangulation all three vertices are of boundary degree $2$ and all three edges have one pre-image in each triangle. In the $(1,1,4)$ triangulation each triangle has a pair of edge pre-images that are identified such that there are two vertices with boundary degree $1$ and the remaining vertex has boundary degree $4$. 

    \begin{figure}
        \centering
        \includegraphics[width=0.5\linewidth]{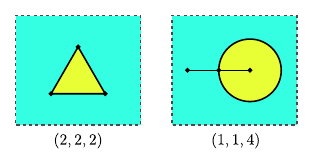}
        \caption{Planar projections of the two triangulations of the $2$-sphere with two triangles, three edges and three vertices.}
        \label{fig:2spheretriangulations}
    \end{figure}
    
    The link of the vertices of boundary degree $1$ is a triangulation of the disc with one 
    boundary edge and vertex. We call this a \emph{monogon} triangulation. The link of the 
    vertices with boundary degree $2$ is a triangulation of the disc with two boundary edges 
    and two vertices. We call this a \emph{bigon} triangulation. The link of the vertex of
    boundary degree $4$ plays a less crucial role in the proof. It follows from the fact 
    that all interior edges of the $3$-ball triangulation must have degree greater than or 
    equal to $3$ that all interior vertices of the link of each boundary vertex must have 
    degree greater than or equal to $3$. 

    First we prove that any monogon triangulation with all interior vertices having degree at 
    least three must have at least three interior vertices. It is sufficient to show that this 
    is true for monogon triangulations that do not have a monogon triangulation as a strict 
    subcomplex. Let $n$ be the number of interior vertices, $F$ the number of triangles,
    and $E$ the number of edges of such a triangulation. Counting edge pre-images yields
    $3F = 2E - 1$ and the Euler characteristic of the disc yields $(n+1) - E + F = 1$. Combining
    these we obtain $F = 2n-1$. Since $F>0$ we deduce $n>0$. If $n=1$ then $F = 1$. The only 
    triangulation of a monogon with a single triangle has an interior vertex of degree $1$ which 
    contradicts our initial assumption. If $n=2$ then $F = 3$. One of these triangles, the 
    \emph{boundary triangle}, must have two of its vertices identified with the boundary vertex. 
    The assumption that there is no monogon strict subcomplex forces the remaining vertex of the 
    boundary triangle to be in the interior of the monogon. This leaves one remaining interior 
    vertex $r$ with all its pre-images in the remaining two triangles. Since $\deg(r) \geq 3$ it 
    follows that one of these triangles has at least two of its vertices identified to $r$ and 
    the edge between them is a loop from $r$ to $r$ in the interior of the monogon. The 
    subcomplex of the triangulation bounded by this loop edge is another monogon triangulation 
    which contradicts our initial assumption. Hence there are at least $3$ interior vertices.

    We prove that any bigon triangulation with all interior vertices having degree at least three must have at least two interior vertices. Let $n, F$ and $E$ be the number of interior vertices, triangles and edges of such a triangulation. Counting edge pre-images yields $3F = 2E -2$ and the Euler characteristic yields $(n+2) - E + F = 1$ which together imply that $2n = F$. There is at least one triangle so we deduce that $n>0$. If $n=1$ then $E = 4$ and $F=2$. In this case the unique interior vertex has degree at least $3$ by assumption. We can exclude the possibility of an edge with both endpoints identified to the interior vertex since such an edge would bound a monogon triangulation that is a strict subcomplex and hence force the existence of other interior vertices. The only remaining possibility is that there are at least three edges with one endpoint at the interior vertex and the other endpoint on either of the boundary vertices. It follows that two of these edges are parallel and bound a strict subcomplex that is another bigon triangulation which again forces the existence of other interior vertices. Hence there are at least $2$ interior vertices.

    We determine the unique bigon triangulation with $2$ interior vertices, each of degree at least $3$. Such a triangulation cannot have either a monogon or another bigon triangulation as a strict subcomplex and so admits no loop edges and no parallel edges except for the two boundary edges. Hence each interior edge $e$ either has one endpoint on the boundary and one endpoint at an interior vertex or both endpoints at distinct interior vertices. For both interior vertices to have a degree of at least $3$, all five combinations of endpoints must appear in the set of interior edges and so the only possible triangulation obtaining the minimum of two interior vertices is the bigon triangulation in Figure~\ref{fig:tetrahedraldisk}. With not much more effort one could prove that the monogon triangulation obtaining the minimum of three interior vertices with each of degree at least $3$ shown in Figure~\ref{fig:tetrahedraldisk} is also unique, but we do not rely on this fact.
    
    \begin{figure}
        \centering
        \includegraphics[width=0.5\linewidth]{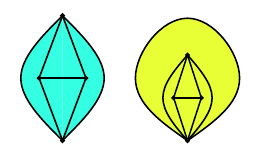}
        \caption{The unique triangulations of a bigon and a monogon minimising interior vertex count for which each interior vertex has degree at least $3$.}
        \label{fig:tetrahedraldisk}
    \end{figure}

    In contrast with the notion of the boundary degree of a vertex, let the \emph{interior degree} $d_{\text{int}}(x)$ of a vertex $x$ be the number of endpoints of interior edges of $H$ identified with $x$. The interior degree of a vertex is equivalent to the number of interior vertices in its link. Equation~\eqref{eq:interiordegree} follows from the fact that each interior edge has two endpoints. 
    
    \begin{equation}\label{eq:interiordegree}
        \sum_{x\in H^{(0)}} d_{\text{int}}(x) = 2e_{\text{int}}^H
    \end{equation}
    
    If $\partial H$ is the $(2,2,2)$ triangulation, all three vertices must have links that are triangulated bigons and so have at least $2$ interior vertices. It follows from Equation~\eqref{eq:interiordegree} that $e_{\text{int}}^H\geq 3$. Since we have $e_{\text{int}}^H<4$ as an assumption this forces $e_{\text{int}}^H=3$ with all three boundary vertices having interior degree $2$. Equivalently, all three vertex links have exactly $2$ interior vertices, and hence are the triangulation of the disc in Figure~\ref{fig:tetrahedraldisk}. 
    
    Now consider a tetrahedron $\mathfrak{t}$ of $H$ with at least one triangle on the boundary. All the vertices of this triangle have links that are the bigon triangulation in Figure~\ref{fig:tetrahedraldisk}. The three edges in $\mathfrak{t}$ incident to the vertex opposite the boundary triangle must all be interior since there exists no triangle in this bigon triangulation with all three vertices on the boundary. It follows that the vertex of $\mathfrak{t}$ opposite the boundary triangle must have a link with a triangle that has three interior vertices. This is not the case for any triangle in the tetrahedral triangulation cut along an edge. Hence no triangulation of a $3$-manifold extending the $(2,2,2)$ triangulation on its boundary exists that has no interior vertices and fewer than four interior edges, each of degree at least three.

    If $\partial H$ is the $(1,1,4)$ triangulation then each of the two boundary degree $1$ vertices must have a vertex link with at least $3$ interior vertices. Again, it follows from Equation~\eqref{eq:interiordegree} that $e_{\text{int}}^H\geq 3$ and so to satisfy our assumption that $e_{\text{int}}^H<4$ there must be exactly three interior edges with every endpoint identified to either of the vertices of boundary degree $1$. The degree $4$ vertex has a vertex link with no interior vertices and so it is a quadrilateral triangulated by the addition of a single diagonal. Two of the vertices of this link are of degree $1$ so these vertices lie on boundary edges of degree $1$. It follows that $d_1^{H}\geq 1$ and hence $e_{\text{int}}^H + d_1^{H} + \frac{d_2^{H}}{2} \geq 4$.

    This proves Inequality~\eqref{eq:equivalentcondition} in the remaining case, where $H$ has a $2$-sphere boundary with two triangles, no interior vertices, and every interior edge has degree at least $3$. This concludes the case analysis verifying Equation~\eqref{eq:generalcase} holds for every component $H$ of $\tilde{\mathcal{C}}$.
    
\end{proof}

\subsection{Summary}
We now summarise the logic completing the proof of Theorem~\ref{thm:main}.

\begin{proof}[Proof of Theorem~\ref{thm:main}]
    The minimal layered triangulations of $L(3n,1)$ defined in \cite{jaco2006layeredtriangulations} provide an upper bound on complexity $c(L(3n,1))\leq 3n-3$ \cite[§2.5, p. 161]{Jaco_2009}. For the lower bound, Proposition~\ref{prop:chargeineqequiv} implies that $c(L(3n,1)) \geq 3n-3$ if and only if Inequality~\eqref{eq:chargeineq} holds. Lemma~\ref{lem:psiineq} implies that Inequality~\eqref{eq:chargeineq} holds if Inequality~\eqref{eq:surplusinequality} (the surplus inequality) holds. We then assume that we have already applied a sequence of $(4-4)$-Pachner moves such that some of the terms in the LHS of Inequality~\eqref{eq:surplusinequality} do not occur. It follows from Lemma~\ref{lem:sufficientconditionpinchedhandle} that the surplus inequality holds whenever Inequality~\eqref{eq:pinchedhandleineq} holds. We then decompose the LHS of Inequality~\eqref{eq:pinchedhandleineq} into the sum of the contribution from each component. The contribution of a component $H$ that is the single tetrahedron layered solid torus $\mathcal{S}(1,2,3)$ is the LHS of Inequality~\eqref{eq:s123case} which directly evaluates to $1$. The contribution of any other component is the LHS of Inequality~\eqref{eq:generalcase} which is non-negative due to Proposition~\ref{prop:desiredbound}. It follows that Inequality~\eqref{eq:pinchedhandleineq} holds globally, yielding the result.
    
\end{proof}

\section{Future Research}
\label{sec:discussion}

The existing literature on applications of the $\mathbb{Z}_2$-Thurston norm to triangulation 
complexity suggests several directions for future research. 

\begin{enumerate}
    \item Unlike~\cite{Jaco_2009} in the case of $L(2n,1)$, our proof does not imply that the minimal layered triangulation of $L(3n,1)$ is the \emph{unique} minimal triangulation. This may be possible with some modification of the above argument that forces the tetrahedra into layers for Inequality~\eqref{eq:chargeineq} to become an equality.
    \item Following~\cite{Rubinstein2024} one could leverage the result from~\cite{hobkirk2026} to derive results on the triangulation complexity of $3$-manifolds with non-trivial second $\mathbb{Z}_3$-homology classes obtained by Dehn filling a link exterior in $S^3$.
    \item We also expect that the same methods used to determine the maximum Euler characteristic of a non-separating $\mathbb{Z}_3$-folded surface embedded in $L(3n,1)$ can be extended to $L(3n,q)$ for all $q$ coprime to $3n$ (with a manageable increase in combinatorial complexity).
    \item Naturally, the methods we have used to determine $c(L(3n,1))$ can then be used to 
    determine more lower bounds on $c(L(3n,q))$. It is possible that these lower bounds will 
    be tight, as was found in~\cite{Jaco_2020} in the case of so-called \textit{balanced lens 
    spaces} using canonical non-orientable surfaces.
    \item Another deeper generalisation is to study $\mathbb{Z}_p$-folded surfaces embedded 
    in $L(pn,q)$ and their relationship with triangulation complexity. The required theory for 
    this has been developed in~\cite{hobkirk2026}, including canonical $\mathbb{Z}_p$-folded 
    surfaces and an appropriate definition of $\chi_p$. In particular, 
    \cite[Conjecture~4.22]{hobkirk2026} and \cite[Corollary~5.3]{hobkirk2026} combined posit 
    that the maximum Euler characteristic of a canonical $\mathbb{Z}_p$-folded surface of a 
    triangulation of $L(pn,1)$ is $2-n$. If this is true then the methods used above may lead 
    to further results on the triangulation complexity of this family of lens spaces. We 
    suspect that bounds on the Euler characteristic of canonical $\mathbb{Z}_p$-folded 
    surfaces may be possible to obtain with methods from link Floer homology.
\end{enumerate}

\bibliographystyle{plain}
\bibliography{refs}

\end{document}